\documentclass[12pt]{amsart}

\usepackage{amsmath}
\usepackage{amssymb}
\usepackage{amsthm}
\usepackage{xcolor}
\usepackage{bm}
\usepackage{graphicx}
\usepackage{pifont}
\usepackage{mathrsfs}
\usepackage{epsfig,verbatim}
\usepackage[pagewise]{lineno}\linenumbers

\newcommand\be{\begin{equation}}
\newcommand\ee{\end{equation}}
\newcommand\bea{\begin{eqnarray}}
\newcommand\eea{\end{eqnarray}}
\newcommand\beaa{\begin{eqnarray*}}
\newcommand\eeaa{\end{eqnarray*}}
\newcommand\beba{\begin{equation}\left\{\begin{array}{rcl}}
\newcommand\eeba{\end{array}\right.\end{equation}}
\newcommand\bebaa{\begin{equation*}\left\{\begin{array}{rcl}}
\newcommand\eebaa{\end{array}\right.\end{equation*}}

\newcommand\R{\mathbb{R}}

\newcommand\bN{{\mathbb{N}}}

\newcommand\cS{{\mathcal{S}}}

\newcommand{\eps}{\epsilon}

\newcommand{\dint}{\displaystyle \int}

\newcommand{\dsup}{\displaystyle \sup}

\newcommand{\dlim}{\displaystyle \lim}
\newcommand{\dmax}{\displaystyle \max}

\newcommand{\dalpha}{D_\xi^\alpha}

\newcommand\ovphi{{\bar{\phi}}}

\newtheorem{theorem}{Theorem}[section]
\newtheorem{lemma}[theorem]{Lemma}
\newtheorem{corollary}[theorem]{Corollary}

\newtheorem{definition}[theorem]{Definition}
\newtheorem{remark}[theorem]{Remark}
\newtheorem{proposition}[theorem]{Proposition}

\numberwithin{equation}{section}

\begin{document}
\nolinenumbers

\title[Asymptotic traveling waves]{Existence and uniqueness of asymptotic traveling waves in a time-fractional bistable reaction-diffusion equation}

\author[A. Kitamachi]{Akihiro Kitamachi}
\address{Department of Mathematics, Hokkaido University, Hokkaido, 060-0810, Japan}
\email{kitamachi.akihiro.k7@elms.hokudai.ac.jp}

\author[H. Ishii]{Hiroshi Ishii}
\address{Research Center of Mathematics for Social Creativity, Research Institute for Electronic Science, Hokkaido University, Hokkaido, 060-0812, Japan}
\email{hiroshi.ishii@es.hokudai.ac.jp}

\thanks{Date: \today. Corresponding author: H. Ishii}

\thanks{{\em Keywords.} {Time-fractional reaction-diffusion equations, Caputo derivative, Traveling wave solutions, Bistable nonlinearity, Continuation method}} 
\thanks{{\em 2020 MSC} Primary {35R11},\ Secondary {35B40, 35C07}}

\begin{abstract}
    We study asymptotic traveling waves for a reaction-diffusion equation with a Caputo time derivative. 
    For general $C^1$ bistable imbalanced nonlinearities, we prove the existence of an asymptotic traveling wave connecting the two stable states. 
    The asymptotic wave speed is unique, and the asymptotic profile is strictly increasing and unique up to translation. 
    The proof uses continuation from a stationary front and the sliding method. 
    For a cubic nonlinearity, we also derive an asymptotic formula for the speed near the balanced case, with an explicit leading coefficient.
    Finally, we perform numerical simulations to compare the speed predicted by the asymptotic formula with the front propagation speed in the original equation.
\end{abstract}

\maketitle
\setcounter{tocdepth}{3}

\section{Introduction}\label{sec:intro}
In this paper, we consider the time-fractional reaction-diffusion equation
\be\label{eq:problem}
    \partial_t^\alpha u = u_{xx} + f(u), \qquad (t>0,\ x\in\R),
\ee
where $u=u(t,x) \in \R$
and $f\in C^1(\R)$ satisfies the following bistable conditions:
there exists $a\in(0,1)$ such that
\be\label{ass:nonlinearity}
\tag{N}
    \begin{gathered}
        f(0)=f(a)=f(1)=0,\quad
        f'(0),f'(1)<0,\quad
        f(u)\ge0\ (u<0),\quad
        f(u)\le0\ (u>1),\\
        f(u)<0\ (u \in (0,a)),\quad
        f(u)>0\ (u\in (a,1)),\quad
        \int_0^1 f(u)\,du>0.
    \end{gathered}
\ee
\noindent
Here, $\partial_t^\alpha$ denotes the Caputo derivative defined by
\beaa
    \partial_t^\alpha u(t,x) 
    := 
    \dfrac{1}{\Gamma(1-\alpha)} 
    \int_0^t 
    \dfrac{\partial u}{\partial s}(s,x)\dfrac{ds}{(t-s)^\alpha},
\eeaa
where $\alpha \in (0,1)$ and $\Gamma$ is the gamma function.

\medskip
In the classical case $\alpha=1$, bistable reaction-diffusion equations arise in models of phase separation~\cite{AC97} and transmission circuits~\cite{NYA65}.
Allen and Cahn~\cite{AC97} modeled interface formation in binary iron-aluminum alloys using the approach developed by Cahn and Hilliard~\cite{CH58}.
This approach is based on free energy and interfacial energy.
Nagumo et al.~\cite{NYA65} proposed an active transmission line with two stable states.
Using mathematical analysis and circuit experiments, they studied the propagation of transition waves between these states and waveform shaping during transmission.
These and related studies have motivated extensive mathematical analysis of scalar bistable reaction-diffusion equations.
For $\alpha=1$, Equation~\eqref{eq:problem} is commonly known as the Allen--Cahn or Nagumo equation.
Such equations are also studied as phase-field models, in which state variables describe the microstructure of materials.

Traveling waves are a characteristic type of solution to reaction-diffusion equations.
They have long been studied to understand phenomena such as propagation and phase separation.
A traveling wave solution in the classical sense has the form $u(t,x)=\phi(x+ct)$ with wave speed $c\in\R$ and profile $\phi$.
Fife and McLeod~\cite{FM77} established the existence, uniqueness, and global stability of front-type traveling wave solutions. 
Research on traveling waves has since expanded to include the construction of entire solutions~\cite{CGN06, FMN04, JM05, MN06}.
In higher dimensions, related work includes the analysis of interface motion~\cite{BBS92, BSS93} and the construction of traveling curved fronts ~\cite{HMR05, HMR06, NT05, NT06}.

Time-fractional reaction-diffusion equations are often used to study spatio-temporal dynamics involving subdiffusion and/or memory effects.
Time-fractional diffusion equations model subdiffusion and can be derived from continuous-time random walks that account for waiting-time effects~\cite{MK00}.
Time-fractional reaction-diffusion systems have also been studied as models of linear reactions under subdiffusive conditions~\cite{HWL06}.
For nonlinear reactions, several time-fractional models have been proposed, including Equation~\eqref{eq:problem} (see, e.g., \cite{AH20, FK26, Nec2010, VNN13}).
In particular, time-fractional bistable reaction-diffusion equations have been studied as phase-field models.
These models incorporate subdiffusion in heterogeneous media and memory effects in materials~\cite{LCWZ18}.
For these equations, the well-posedness of initial-boundary value problems, maximum principles, and numerical approximations have been studied~\cite{DYZ20}.
An energy dissipation law of integral type has also been established~\cite{TYZ19}.

To investigate the effects of fractional time derivatives on propagation, we consider traveling waves.
However, the usual traveling wave reduction does not directly apply to Equation~\eqref{eq:problem} with the Caputo derivative.
Indeed, substituting $u(t,x)=\phi(x+ct)$ into the definition of the Caputo derivative gives
\beaa
    \partial^{\alpha}_{t} u(t,x) = \dfrac{c}{\Gamma(1-\alpha)} \int^{t}_{0} \dfrac{\phi'(x+cs)}{(t-s)^{\alpha}}\, ds= \dfrac{c}{\Gamma(1-\alpha)} \int^{t}_{0} \dfrac{\phi'(x+c(t-s))}{s^{\alpha}}\, ds.
\eeaa
Introducing the moving coordinate $\xi:=x+ct$, we obtain
\be\label{eq:Caputo}
\partial^{\alpha}_{t} u(t,x)= 
\begin{cases}
    \dfrac{c^{\alpha}}{\Gamma(1-\alpha)} \dint^{ct}_{0} \dfrac{\phi'(\xi-s)}{s^{\alpha}}\, ds &(c>0), \\
    -\dfrac{|c|^{\alpha}}{\Gamma(1-\alpha)} \dint^{|c|t}_{0} \dfrac{\phi'(\xi+s)}{s^{\alpha}}\,ds &(c<0).
\end{cases}
\ee
Thus, the Caputo derivative retains an explicit dependence on $t$ in the moving coordinate.
This prevents the usual reduction to a stationary profile equation.

This motivates the formal derivation of an equation for a possible asymptotic profile. 
Specifically, we assume that
\beaa
     u(t,x) \simeq \phi(x+ct) \quad (t\to+\infty).
\eeaa
Keeping $\xi$ fixed, we then formally let $t\to+\infty$ on the right-hand side of \eqref{eq:Caputo} to obtain
\beaa
\partial^{\alpha}_{t} u(t,x) \simeq 
\begin{cases}
    \dfrac{c^{\alpha}}{\Gamma(1-\alpha)} \dint^{+\infty}_{0} \dfrac{\phi'(\xi-s)}{s^{\alpha}}\, ds &(c>0), \\
    -\dfrac{|c|^{\alpha}}{\Gamma(1-\alpha)} \dint^{+\infty}_{0} \dfrac{\phi'(\xi+s)}{s^{\alpha}}\,ds &(c<0).
\end{cases}
\eeaa
This formal procedure removes the explicit dependence on $t$ and yields a profile equation.
A pair $(c,\phi)$ satisfying this profile equation is called an asymptotic traveling wave~\cite{Ishii2025}.
Asymptotic traveling waves have been analyzed for time-fractional reaction-diffusion equations~\cite{Ishii2025, Nec2010, VNN13}.
Such waves have also been studied for time-fractional neural field models~\cite{Gonzalez, Ishii-preprint}.
In particular, for a time-fractional neural field model, the second author has obtained global-in-time solutions that converge to asymptotic traveling waves~\cite{Ishii-preprint}.
These studies examine the effects of fractional time derivatives by analyzing wave speeds and profiles.

The asymptotic profile equation is nonlocal in the spatial variable.
Its analysis is therefore also related to traveling wave theory for spatially nonlocal reaction-diffusion equations.
For spatially nonlocal equations with bistable reaction terms, the existence and uniqueness of traveling wave solutions have been studied (see, e.g., \cite{AK15, Bates1997, Chen97, GUI15}).
The existence and uniqueness of traveling wave solutions have also been established for reaction-diffusion equations with memory terms~\cite{MR24}.

Asymptotic traveling waves for Equation~\eqref{eq:problem} were studied in \cite{Nec2010} for the nonlinearities
\be\label{Nec1}
    f(u) = 
    \begin{cases}
        -ku &(u<a), \\
        k(1- u) &(u> a),
    \end{cases}
    \qquad (0<a<1,\ k>0)
\ee
and
\be\label{Nec2}
    f(u) = u - u^3 + \mu\qquad (|\mu|\ll 1).
\ee
For \eqref{Nec1}, several properties of the wave speed and profile were obtained using a solution representation.
Such properties were also derived for \eqref{Nec2} through a perturbative argument in $\mu$.
However, the previous study \cite{Nec2010} does not provide a theorem establishing the existence, monotonicity, and uniqueness of fronts for general bistable nonlinearities.
It also does not examine the relationship between asymptotic traveling waves and the long-time behavior of solutions to Equation~\eqref{eq:problem}.

In this paper, we analyze the asymptotic profile equation derived from Equation~\eqref{eq:problem} for general $C^1$ bistable nonlinearities satisfying \eqref{ass:nonlinearity}.
For each $\alpha\in(0,1)$, we establish the existence of an asymptotic traveling wave $(c,\phi)$ with $c>0$ whose profile connects the two stable equilibria.
This result does not require $\int_0^1 f(u)\,du$ to be small.
The wave speed is unique, and the profile is unique up to translation.
We also prove that every asymptotic traveling wave profile connecting $0$ at $-\infty$ to $1$ at $+\infty$ is strictly increasing.
Furthermore, we use numerical simulations of Equation~\eqref{eq:problem} to examine the relationship between asymptotic traveling waves and the long-time behavior of solutions.

The paper is organized as follows. 
Section~\ref{sec:main} states our main results on the existence and uniqueness of asymptotic traveling waves.
Section~\ref{sec:pre} establishes the preliminary lemmas needed for the proofs.
Section~\ref{sec:mono-uni} proves the uniqueness of asymptotic traveling waves and the strict monotonicity of their profiles. 
Section~\ref{sec:proof} shows the existence of asymptotic traveling waves. 
Section~\ref{sec:num} derives an asymptotic formula for the wave speed in the case of a cubic nonlinearity and presents numerical simulations to investigate the long-time behavior of solutions to Equation~\eqref{eq:problem} and its relationship to asymptotic traveling waves.
Finally, Section~\ref{sec:dis} discusses remaining open problems.

\bigskip
\section{Main result}\label{sec:main}

We study asymptotic traveling waves connecting the two stable equilibria $0$ and $1$. 
Under the imbalance condition $\int_0^1 f(u)\,du>0$, it suffices to consider asymptotic traveling waves with positive speed $c$ (see Remark~\ref{rem:speed}). 
Accordingly, we adopt the following definition:
\begin{definition}
    A pair $(c,\phi)\in (0,+\infty) \times C^2(\R)$ is an asymptotic traveling wave of Equation~\eqref{eq:problem} if it satisfies 
    \be
    \label{eq:atw}
        \frac{c^\alpha}{\Gamma(1-\alpha)} 
        \int_{0}^{+\infty}
        \frac{\phi'(\xi-s)}{s^\alpha}\,ds
        =
        \phi''(\xi)+f(\phi(\xi))
        \qquad (\xi\in\R).
    \ee
    We call $c$ the asymptotic wave speed and $\phi$ the asymptotic profile.
\end{definition}

\medskip
We focus on asymptotic traveling waves satisfying
\be
    \label{ass:boundary_condition}
    \tag{B}
    \dlim_{\xi\to-\infty}\phi(\xi)=0,
    \qquad
    \dlim_{\xi\to+\infty}\phi(\xi)=1.
\ee
We then obtain the following theorem.

\begin{theorem}\label{thm:main}
There exists an asymptotic traveling wave $(c,\phi)\in (0,+\infty)\times C_b^2(\R)$ of Equation~\eqref{eq:problem} satisfying \eqref{ass:boundary_condition}.
The speed $c$ is unique, and the profile $\phi$ is unique up to translation.
Moreover, $\phi'>0$.
\end{theorem}

\medskip
The positivity of $c$ and $\phi'$ yields the following sub-solution property, as in~\cite[Corollary 2.3]{Ishii2025}.
\begin{corollary}\label{cor:sub}
    For $\alpha\in(0,1)$, let $(c,\phi)$ be an asymptotic traveling wave given by Theorem~\ref{thm:main}.
    Then $v(t,x)=\phi(x+ct)$ is a sub-solution of Equation~\eqref{eq:problem}. 
    More precisely,
    \beaa
    \partial^{\alpha}_t v < v_{xx} + f(v) \quad (t>0,\ x\in\R).
    \eeaa
\end{corollary}

\medskip
Corollary~\ref{cor:sub} shows that asymptotic traveling waves can be used as comparison functions for the original evolution problem.
Indeed, if a bounded solution $u$ of Equation~\eqref{eq:problem} satisfies
\beaa
    u(0,x)\ge \phi(x+h)\qquad(x\in\R)
\eeaa
for some $h\in\R$, then the comparison principle yields
\beaa
    u(t,x)\ge \phi(x+ct+h)\qquad(t>0,\ x\in\R).
\eeaa
Thus, $\phi(x+ct+h)$ gives a lower bound for solutions whose initial data satisfy $u(0,x)\ge\phi(x+h)$ for all $x\in\R$.
In this sense, the analysis of asymptotic traveling waves provides insight into the propagation behavior of solutions to Equation~\eqref{eq:problem}.
We conjecture that, for suitable initial data, global-in-time
solutions to Equation~\eqref{eq:problem} converge to asymptotic traveling waves, but this remains an open problem.

Theorem~\ref{thm:main} establishes existence, monotonicity, and uniqueness for asymptotic traveling waves of time-fractional reaction-diffusion equations.
These results are analogous to those known for classical bistable traveling waves.
For nonlinearities with a small imbalance, we can also examine how the asymptotic wave speed depends on the fractional order $\alpha$ (see Remark~\ref{rem:speed-zero}).

In Section~\ref{sec:num}, we numerically investigate front propagation in Equation~\eqref{eq:problem}.
For a cubic nonlinearity, we derive an asymptotic formula for the speed of asymptotic traveling waves near the balanced case.
This formula describes how the fractional order enters the speed scaling (see Section~\ref{subsec:sp-atw}).
For the cubic case, Nec et al.~\cite{Nec2010} derived the $1/\alpha$ power law for the wave speed near the balanced case.
Their perturbative calculation in a small external field also yielded an integral representation of the speed coefficient.
Here, we derive the corresponding asymptotic formula from solutions constructed using the implicit function theorem.
We express the speed coefficient explicitly in terms of the gamma function and the Riemann zeta function.
We then compare the numerically computed front speeds with those predicted by the asymptotic formula. 
We use this comparison to examine how the fractional order and the degree of imbalance affect front propagation.

In proving Theorem~\ref{thm:main}, we cannot directly use classical phase-plane methods because Equation~\eqref{eq:atw} is nonlocal.
Although traveling waves for nonlocal equations have been widely studied, most methods have been developed for equations with ordinary first-order time derivatives.
To the best of our knowledge, no existing result directly establishes existence, monotonicity, and uniqueness for Equation~\eqref{eq:atw}.
To prove existence, we use a continuation method, as in~\cite{Bates1997, EM93}.
The starting point is a stationary front for a bistable nonlinearity $g$ satisfying $\int_0^1 g(u)\,du=0$.
The nonlinearity is then continuously deformed from $g$ to $f$.
In a framework based on $C_0(\R)$, we construct solutions locally in the continuation parameter using the implicit function theorem. 
We obtain solutions at limiting parameter values by passing to the limit along sequences of solutions.
We prove monotonicity and uniqueness using the sliding method, following the ideas in~\cite{BN91, GUI15}.

For notational convenience, we set
\beaa
    \dalpha \psi(\xi)
    :=
    \frac{1}{\Gamma(1-\alpha)}
    \int_0^{+\infty} \frac{\psi'(\xi-s)}{s^\alpha}\, ds,
    \qquad
    \theta := c^\alpha.
\eeaa
Then Equation~\eqref{eq:atw} is rewritten as
\be\label{eq:atw-theta}
    \theta\dalpha \phi(\xi) = \phi''(\xi)+f(\phi(\xi))
    \qquad (\xi\in\R).
\ee

\bigskip
\section{Preliminaries}\label{sec:pre}
We present some preliminary results used in the proof of the main theorem.
In Section~\ref{subsec:stasol}, we recall the existence of a stationary front in the balanced case and state the solvability condition for the equation linearized about this front.
In Section~\ref{subsec:frac}, we collect basic lemmas concerning the fractional derivative $D_\xi^\alpha$ and the profile equation, including maximum principles, regularity results, and lemmas for passing to the limit.

For each $k\in\bN \cup \{0\}$, we define the spaces
\beaa
    C_b^k(\R)
    &:=&
    \left\{v\in C^k(\R):\,v^{(j)} \in L^\infty(\R) \quad (0 \le j \le k) \right\},\\
    C_0^k(\R)
    &:=&
    \{v\in C_b^k(\R):\,\dlim_{\xi \to \pm \infty} v^{(j)}(\xi) = 0 \quad (0 \le j \le k)\}.
\eeaa
We equip these spaces with the norm
\beaa
    \|v\|_{C_b^k} := \sum_{j=0}^k\|v^{(j)}\|_{L^\infty}.
\eeaa
Also, we write $C_b^0(\R)=C_b(\R)$ and $C_0^0(\R)=C_0(\R)$.
Note that both $C_b^k(\R)$ and $C_0^k(\R)$ are Banach spaces equipped with this norm.
\medskip
\subsection{Stationary fronts in the balanced case}\label{subsec:stasol}
In this subsection, we consider the stationary problem
\be\label{eq:balanced}
    \phi_0''+g(\phi_0)=0.
\ee
Let $a\in(0,1)$ be as in assumption~\eqref{ass:nonlinearity}.
We assume that $g\in C^1(\R)$ satisfies the following balanced bistable conditions:
\beaa
    \begin{aligned}
    g(0)=g(a)=g(1)=0,\qquad
    g'(0)<0,&\qquad
    g'(1)<0,\\
    g(u)<0 \quad (u\in(0,a)),\qquad
    g(u)>0 \quad (u\in(a,1)),&\qquad
    \int_0^1 g(u)\,du = 0.
    \end{aligned}
\eeaa
Such a function $g$ can be constructed directly from the nonlinearity $f$ satisfying \eqref{ass:nonlinearity} (see Section~\ref{sec:proof}).

In this subsection, we discuss the properties of stationary fronts of Equation~\eqref{eq:balanced} satisfying \eqref{ass:boundary_condition} and of the associated linearized operator.
We first recall the following known result on the existence of a stationary front.
\begin{proposition}\cite[Theorem~3.9]{Volpert1994}\label{prop:stationary-front}
    Equation~\eqref{eq:balanced} has a solution $\phi_0\in C_b^2(\R)$ satisfying \eqref{ass:boundary_condition} and $\phi_0'>0$.
\end{proposition}

Note that $\phi_0'$ is uniformly continuous since $\phi_0'\in C_b^1(\R)$.
Moreover, it follows from \eqref{ass:boundary_condition} and the positivity of $\phi_0'$ that $\phi_0'$ is integrable on $\R$.
Hence $\phi_0'\in C_0(\R)$.
Furthermore, since $g\in C^1(\R)$, Equation~\eqref{eq:balanced} yields
\beaa
    \phi_0''=-g(\phi_0),\qquad \phi_0'''=-g'(\phi_0)\phi_0'.
\eeaa
Thus both $\phi_0''$ and $\phi_0'''$ tend to zero as $\xi\to\pm\infty$, and consequently $\phi_0'\in C_0^2(\R)$.

We define the linearized operator $L_0:C_0^2(\R)\to C_0(\R)$  about the stationary front $\phi_0$ by
\beaa
    L_0 h:=h''+g'(\phi_0)h\qquad (h\in C_0^2(\R)).
\eeaa
Differentiating Equation~\eqref{eq:balanced} with respect to $\xi$, we obtain $L_0\phi_0'=0$.
The following lemma gives a solvability condition in $C_0^2(\R)$ for the linearized equation with a normalization condition.

\begin{lemma}\label{lem:solvability-L0}
    For any $(A,b)\in C_0(\R)\times\R$, the problem
    \be
        \label{eq:nonhomo}
        L_0h = A, \qquad h(0) = b,
    \ee
    has a solution $h\in C_0^2(\R)$ if and only if
    \be\label{condi:solve}
        \int_{\R}A(\xi)\phi_0'(\xi)\,d\xi=0.
    \ee
    Moreover, the solution, if it exists, is unique.
\end{lemma}

\begin{proof}
    We follow the argument in \cite[Lemma~2.2]{Alfaro2008}.
    Fix an arbitrary $(A,b)\in C_0(\R)\times\R$.

    First, suppose that \eqref{eq:nonhomo} has a solution $h\in C_0^2(\R)$.
    Multiplying $L_0h=A$ by $\phi_0'$ and integrating by parts over $\R$, we obtain
    \beaa
        \int_{\R}A(\xi)\phi_0'(\xi)\,d\xi
        =\int_{\R}L_0h(\xi)\phi_0'(\xi)\,d\xi
        =\int_{\R}h(\xi)L_0\phi_0'(\xi)\,d\xi
        =0.
    \eeaa
    Hence \eqref{condi:solve} holds.

    Conversely, assume that \eqref{condi:solve} holds.
    Using the representation formula in \cite[formula~(2.12)]{Alfaro2008}, we define
    \beaa
        \begin{aligned}
            h(\xi)
            &=\phi_0'(\xi)\left[
                \frac{b}{\phi_0'(0)}
                +\int_0^\xi\frac{J(s)}{(\phi_0'(s))^2}\,ds
            \right],\\
            \quad
            J(\xi)
            &:=\int_{-\infty}^{\xi}A(s)\phi_0'(s)\,ds
            =-\int_{\xi}^{\infty}A(s)\phi_0'(s)\,ds.
        \end{aligned}
    \eeaa
    Since $\phi_0'>0$, $\phi_0'\in C^2(\R)$, and $A\in C_0(\R)$, this formula defines a function $h\in C^2(\R)$.
    By construction, $h(0)=b$, and direct differentiation gives
    \be\label{eq:nonhomo-pre}
        h''+g'(\phi_0)h=A.
    \ee
    It remains to show that $h\in C_0^2(\R)$.
    For every $\xi\in\R$, we have
    \beaa
        |J(\xi)|
        \le
        \sup_{s\ge\xi}|A(s)|
        \bigl(1-\phi_0(\xi)\bigr),
        \qquad
        |J(\xi)|
        \le
        \sup_{s\le\xi}|A(s)|\,\phi_0(\xi).
    \eeaa
    Since $g(0)=g(1)=0$ and $g$ is differentiable at $0$ and $1$, we also have
    \beaa
        \lim_{\xi\to+\infty}
        \frac{g(\phi_0(\xi))}{1-\phi_0(\xi)}
        =-g'(1)>0,\qquad
        \lim_{\xi\to-\infty}
        \frac{g(\phi_0(\xi))}{\phi_0(\xi)}
        =g'(0)<0.
    \eeaa
    Together with $A\in C_0(\R)$, these estimates imply that
    $J(\xi)/g(\phi_0(\xi))\to0$ as $\xi\to\pm\infty$.
    Since $1/\phi_0'(\xi)\to+\infty$ as $\xi\to\pm\infty$, l'H\^opital's rule yields
    \beaa
        \begin{aligned}
        \lim_{\xi\to\pm\infty}h(\xi)
        =
        \lim_{\xi\to\pm\infty}
        \frac{\dfrac{b}{\phi'_0(0)}+\dint_0^\xi\frac{J(s)}{(\phi'_0(s))^2}\,ds}{1/\phi'_0(\xi)}
        =
        \lim_{\xi\to\pm\infty} \frac{-J(\xi)}{\phi''_0(\xi)}
        =
        \lim_{\xi\to\pm\infty} \frac{J(\xi)}{g(\phi_0(\xi))}
        =0.
        \end{aligned}
    \eeaa
    Thus $h\in C_0(\R)$.
    Equation~\eqref{eq:nonhomo-pre} and the boundedness of $g'(\phi_0)$ then imply that $h''\in C_0(\R)$.
    Furthermore, Taylor's theorem gives
    \beaa
        h'(\xi)
        =h(\xi+1)-h(\xi)
        -\int_0^1(1-t)h''(\xi+t)\,dt.
    \eeaa
    Since both $h$ and $h''$ tend to zero as $\xi\to\pm\infty$, this identity implies that $h'(\xi)\to0$ as $\xi\to\pm\infty$.
    Hence $h\in C_0^2(\R)$.

    Finally, we prove uniqueness.
    Let $h_1,h_2\in C_0^2(\R)$ satisfy \eqref{eq:nonhomo}, and set $w:=h_1-h_2$.
    Since $L_0w=L_0\phi_0'=0$, we have
    \beaa
        \bigl(\phi_0'w'-\phi_0''w\bigr)'=0.
    \eeaa
    Thus $\phi_0'w'-\phi_0''w$ is constant, and its decay at infinity implies that it is zero.
    Since $\phi_0'>0$, the quotient $w/\phi_0'$ is constant.
    Together with $w(0)=0$, this gives $w\equiv0$.
    Therefore, the solution is unique.
\end{proof}

\medskip
Lemma~\ref{lem:solvability-L0} gives the solvability condition for the linearized equation. 
When this condition is satisfied, the normalization $h(0)=b$ determines the solution uniquely. 
In Section~\ref{sec:proof}, we treat $\theta=c^\alpha$ as an additional unknown and include the normalization condition to form an augmented linearized operator. 
Using Lemma~\ref{lem:solvability-L0}, we prove that this operator is invertible. 
We then apply the implicit function theorem to construct asymptotic traveling waves near the stationary front in the balanced case.

\medskip
\subsection{Properties of the fractional derivative}\label{subsec:frac}

In this subsection, we collect some basic properties of the fractional derivative $\dalpha$.
We first establish a maximum principle using the sign of the fractional derivative at a global minimum point.
We then discuss the boundedness of $\dalpha$, its commutation with differentiation, and the decay of fractional derivatives at infinity.
Finally, we discuss the interchange of fractional differentiation and limits under local convergence.

We first present a lemma on the sign of the fractional derivative at a global minimum point.
This lemma provides the basis for the maximum principle.
\begin{lemma}\cite[Lemma~3.6]{Ishii2025}\label{lem:min-sgn}
    Let $\psi \in C_b^1(\R)$.
    Assume that $\psi(\xi)$ attains its global minimum at $\xi^* \in \R$, that is, $\psi(\xi^*) = \min_{\xi \in \R}\psi(\xi)$.
    Then, we have
    \beaa
        \dalpha \psi(\xi^*)
        =
        -\dfrac{\alpha}{\Gamma(1-\alpha)} 
        \int_{0}^{+\infty} 
        \dfrac{\psi(\xi^*-s) - \psi(\xi^*)}{s^{1+\alpha}} \, ds
        \le 
        0.
    \eeaa
    Moreover, $\dalpha \psi(\xi^*) =0$
    if and only if $\psi(\xi) = \psi(\xi^*)$ for any $\xi \le \xi^*$.
\end{lemma}

\medskip
Although Lemma~\ref{lem:min-sgn} is stated for a minimum point, the corresponding result at a maximum point follows by applying it to $-\psi$.
Using this sign property of $\dalpha$, we obtain the following maximum principle.
\begin{lemma}\label{lem:max}
    Let $d \in L^\infty(\R)$ and $\theta >0$.
    Assume that $w\in C_b^2(\R)$ satisfies
    \beaa
            w''(\xi) - \theta \dalpha w(\xi) + d(\xi)w(\xi) \le 0 \qquad (\xi \in \R),
            \qquad
            \displaystyle\liminf_{\xi\to \pm \infty}  w(\xi) \ge 0.
    \eeaa
    If there exists a set $H \subset \R$ such that
    $w(\xi) \ge 0$ for $\xi \in H$ and 
    $d(\xi) \le 0$ for $\xi \not \in H$, 
    then either $w>0$ or $w \equiv 0$.
\end{lemma}

\begin{proof}
    First, we show that $w\ge0$.
    The conclusion is immediate if $w\equiv0$, so we assume that $w\not\equiv0$.
    Suppose, for contradiction, that $w$ takes a negative value.
    By continuity and the condition at infinity, $w$ attains a negative global minimum at some $\xi_1\in\R$, that is,
    $w(\xi_1)=\min_{\xi\in\R}w(\xi)<0$.
    Then $w''(\xi_1)\ge0$, and since $\xi_1\notin H$, we have $d(\xi_1)\le0$.
    Moreover, Lemma~\ref{lem:min-sgn} gives $\dalpha w(\xi_1)\le0$.
    If equality held, we would have $w(\xi)=w(\xi_1)<0$ for all $\xi\le\xi_1$, contradicting the condition at infinity.
    Thus $\dalpha w(\xi_1)<0$, and hence
    \beaa
    w''(\xi_1)-\theta\dalpha w(\xi_1)
    +d(\xi_1)w(\xi_1)>0.
    \eeaa
    This contradicts the assumed differential inequality.
    Therefore, $w\ge0$.

    Next, we show that if $w(\xi_0)=0$ for some $\xi_0\in\R$, then $w\equiv0$.
    Since $\xi_0$ is a global minimum point of $w$, we have $w'(\xi_0)=0$ and $w''(\xi_0)\ge0$.
    Lemma~\ref{lem:min-sgn} also gives $\dalpha w(\xi_0)\le0$.
    Consequently,
    \beaa
        0\ge
        w''(\xi_0)-\theta\dalpha w(\xi_0)
        \ge-\theta\dalpha w(\xi_0)\ge0,
    \eeaa
    which implies that $\dalpha w(\xi_0)=0$.
    Applying Lemma~\ref{lem:min-sgn} again, we obtain $w(\xi)=0$ for all $\xi\le\xi_0$.
    
    It remains to show that $w$ also vanishes for $\xi>\xi_0$.
    Since $w(\xi_0)=w'(\xi_0)=0$, we obtain
    \be\label{ineq:w}
        0\le w(\xi)
        =
        \int_{\xi_0}^{\xi}(\xi-\tau)w''(\tau)\,d\tau
        \le
        \theta\int_{\xi_0}^{\xi}(\xi-\tau)\dalpha w(\tau)\,d\tau
        +\|d\|_{L^\infty}\int_{\xi_0}^{\xi}(\xi-\tau)w(\tau)\,d\tau.
    \ee
    For the first integral on the right-hand side, we use the fact that $w(\xi)=0$ for $\xi\le\xi_0$.
    By interchanging the order of integration and integrating by parts, we obtain
    \be\label{integ:frac}
        \int_{\xi_0}^\xi(\xi-\tau)\dalpha w(\tau)\,d\tau
        =
        \frac{1}{\Gamma(2-\alpha)}
        \int_{\xi_0}^\xi(\xi-s)^{1-\alpha}w(s)\,ds.
    \ee
    For $\ell>0$, set $M:=\dmax_{\xi\in[\xi_0,\xi_0+\ell]}w(\xi)$.
    Combining \eqref{ineq:w} and \eqref{integ:frac}, we obtain
    \beaa
        0\le w(\xi)
        \le
        \left[
            \frac{\theta(\xi-\xi_0)^{2-\alpha}}
                 {\Gamma(3-\alpha)}
            +\frac{\|d\|_{L^\infty}}{2}(\xi-\xi_0)^2
        \right]M
    \eeaa
    for every $\xi\in[\xi_0,\xi_0+\ell]$.
    Taking the maximum over this interval gives
    \beaa
        0\le M\le
        \left[
            \frac{\theta\ell^{2-\alpha}}{\Gamma(3-\alpha)}
            +\frac{\|d\|_{L^\infty}}{2}\ell^2
        \right]M.
    \eeaa
    Choosing $\ell>0$ sufficiently small that the expression in brackets is less than $1$, we obtain $M=0$.
    Since this choice of $\ell$ is independent of $\xi_0$, we can repeat the argument on successive intervals of length $\ell$ to conclude that $w\equiv0$.
\end{proof}

\medskip
We next establish the regularity and decay properties needed to differentiate the profile equation and pass to the limit.
\begin{lemma}\label{lem:D-bdd}
    For all $k\in \bN$, the following statements hold:
    \begin{enumerate}
        \item[(i)]
        for all $v\in C_b^k(\R)$, 
        $\dalpha v\in C_b^{k-1}(\R)$, $\|\dalpha v\|_{C_b^{k-1}}\le \dfrac{1}{\Gamma(1-\alpha)} \left( \dfrac{1}{1-\alpha} + 2 \right)\|v\|_{C_b^k}$.

        \item[(ii)]
        for all $v\in C_b^k(\R)$,
        $(\dalpha v)^{(j)}=\dalpha(v^{(j)})\,\,\,(j=0,\ldots,k-1)$.

        \item[(iii)]
        if $v\in C_0^k(\R)$, then $\dalpha v\in C_0^{k-1}(\R)$.
    \end{enumerate}
\end{lemma}

\begin{proof}
    We first derive an estimate for $k=1$.
    We then use this estimate and differentiation under the integral sign to prove (i) and (ii).
    Let $v\in C_b^1(\R)$.
    We split the integral into the parts over $(0,1)$ and $(1,\infty)$.
    Integrating the latter part by parts and using the boundedness of $v$, we obtain
    \be\label{eq:Dalpha-split}
        \Gamma(1-\alpha)\dalpha v(\xi)
        =
        \int_0^1 s^{-\alpha}v'(\xi-s)\,ds
        +v(\xi-1)
        -\alpha\int_1^\infty s^{-1-\alpha}v(\xi-s)\,ds.
    \ee
    Since $s^{-\alpha}\in L^1(0,1)$ and $s^{-1-\alpha}\in L^1(1,\infty)$,
    the dominated convergence theorem implies that $\dalpha v$ is continuous.
    Moreover,
    \be\label{ineq:Dalpha-est}
        \|\dalpha v\|_{C_b}
        \le
        \frac{1}{\Gamma(1-\alpha)}
        \left(
        \frac{\|v'\|_{C_b}}{1-\alpha} + 2\|v\|_{C_b}
        \right).
    \ee
    This proves (i) and (ii) for $k=1$.

    Next, let $k\ge2$ and $v\in C_b^k(\R)$.
    The derivatives of $v$ up to order $k$ are bounded, and the kernels above are integrable.
    Thus, we can differentiate the right-hand side of \eqref{eq:Dalpha-split} under the integral sign up to $k-1$ times.
    The resulting expressions agree with those obtained by applying \eqref{eq:Dalpha-split} to $v^{(j)}$.
    Hence
    \beaa
        (\dalpha v)^{(j)}=\dalpha(v^{(j)})
        \qquad (j=0,\ldots,k-1),
    \eeaa
    which proves (ii).
    Applying the result for $k=1$ to each $v^{(j)}\in C_b^1(\R)$, we obtain $\dalpha v\in C_b^{k-1}(\R)$.
    The estimate \eqref{ineq:Dalpha-est} then gives (i).
    
    Finally, we prove (iii).
    Let $v\in C_0^k(\R)$.
    Then $v^{(j)},v^{(j+1)}\in C_0(\R)$ for each $j=0,\ldots,k-1$.
    Applying \eqref{eq:Dalpha-split} to $v^{(j)}$ and using the dominated convergence theorem together with (ii), we obtain
    \beaa
        (\dalpha v)^{(j)}(\xi)
        =\dalpha(v^{(j)})(\xi)\to0
        \qquad (\xi\to\pm\infty).
    \eeaa
    Hence $\dalpha v\in C_0^{k-1}(\R)$.
\end{proof}

\medskip
We next examine the behavior at infinity of the derivative and the fractional derivative of a monotone function.
\begin{lemma}\label{lem:Dalpha-phi-C0}
    Let $v\in C_b^2(\R)$. If $v$ is monotone, then $\dalpha v\in C_0(\R)$.
\end{lemma}
\begin{proof}
    It suffices to consider the case $v'\ge0$.
    Since $v$ is bounded and nondecreasing, the limits $L_\pm:=\lim_{\xi\to\pm\infty}v(\xi)$ exist.
    Hence $v'\in L^1(\R)$.
    Moreover, $v\in C_b^2(\R)$ implies that $v'$ is uniformly continuous on $\R$.
    Together with its integrability, this yields $v'\in C_0(\R)$.
    
    By Lemma~\ref{lem:D-bdd}, $\dalpha v$ is continuous.
    Moreover, \eqref{eq:Dalpha-split} and the dominated convergence theorem imply that, as $\xi\to\pm\infty$,
    \beaa
        \Gamma(1-\alpha)\dalpha v(\xi)
        \to
        0+L_\pm - \alpha L_\pm\int_1^{+\infty} s^{-1-\alpha}\,ds
        =
        0.
    \eeaa
    Hence $\dalpha v\in C_0(\R)$.
\end{proof}

\medskip
Finally, we prove a lemma on passing to the limit in $\dalpha$ under local convergence.
This will be used to show that the limit of a sequence of solutions satisfies the profile equation.
\begin{lemma}\label{lem:Dalpha-local-conv}
    Let $v_n,v\in C_b^1(\R)$ and assume that
    $\dsup_{n\in\bN}\|v_n\|_{C_b}<\infty$.
    Then
    \beaa
        v_n\to v
        \quad\text{in }C_{\mathrm{loc}}^1(\R)
        \quad\Longrightarrow\quad
        \dalpha v_n\to\dalpha v
        \quad\text{in }C_{\mathrm{loc}}(\R).
    \eeaa
\end{lemma}

\begin{proof}
    Set $w_n:=v_n-v$.
    Then $w_n\to0$ in $C_{\mathrm{loc}}^1(\R)$ and $M:=\sup_{n\in\bN}\|w_n\|_{C_b}<\infty$.
    Fix an arbitrary compact set $K\subset\R$, and let $R>1$.
    Set $K_R:=K-[0,R]$.
    We apply \eqref{eq:Dalpha-split} to $w_n$.
    Splitting the last integral at $s=R$, we obtain
    \beaa
        \Gamma(1-\alpha)
        \sup_{\xi\in K}|\dalpha w_n(\xi)|
        \le
        \frac{1}{1-\alpha}\sup_{\xi\in K_R}|w_n'(\xi)|
        +2\sup_{\xi\in K_R}|w_n(\xi)|
        +MR^{-\alpha}.
    \eeaa
    Since $K_R$ is compact for each fixed $R$, letting $n\to\infty$ yields
    \beaa
        \limsup_{n\to\infty}
        \sup_{\xi\in K}|\dalpha w_n(\xi)|
        \le
        \frac{M}{\Gamma(1-\alpha)}R^{-\alpha}.
    \eeaa
    Letting $R\to\infty$, we conclude that
    $\sup_{\xi\in K}|\dalpha w_n(\xi)|\to0$ as $n\to\infty$.
    Since $K$ was arbitrary,
    $\dalpha v_n\to\dalpha v$ in $C_{\mathrm{loc}}(\R)$.
\end{proof}

\bigskip
\section{Monotonicity and uniqueness of asymptotic traveling waves}\label{sec:mono-uni}

In this section, we establish some basic properties of asymptotic traveling waves.
We first show that it suffices to consider the case $\theta>0$ when studying Equation~\eqref{eq:atw-theta}.
\begin{proposition}\label{prop:posi-speed}
    If $(\theta,\phi)\in\R\times C_b^2(\R)$ is a solution of Equation~\eqref{eq:atw-theta} satisfying \eqref{ass:boundary_condition}, then $\theta>0$.
\end{proposition}
\begin{proof}
    Let $(\theta,\phi)\in\R\times C_b^2(\R)$ be a solution of Equation~\eqref{eq:atw-theta} satisfying \eqref{ass:boundary_condition}.
    Multiplying the equation by $\phi'$, we obtain
    \be\label{eq:pre-energy}
    \theta \phi'(\xi) \dalpha \phi(\xi) = \dfrac{d}{d\xi} \left[ \dfrac{1}{2} (\phi'(\xi))^2 + F(\phi(\xi)) \right],
    \ee
    where
    \beaa
    F(u) := \int^{u}_{0} f(s)\,ds.
    \eeaa
    Since $\phi\in C_b^2(\R)$ satisfies \eqref{ass:boundary_condition}, we have $\phi'(\pm\infty)=0$.

    We now rewrite the left-hand side of \eqref{eq:pre-energy}.
    This calculation follows the argument for the Caputo derivative $\partial_t^\alpha$ in~\cite[Lemma~3.1]{QTWY}.
    We include the derivation for $\dalpha$.
    First, integration by parts gives
    \beaa
        \dalpha \phi(\xi) = \dfrac{\alpha}{\Gamma(1-\alpha)} \int^{+\infty}_{0} \dfrac{\phi(\xi)-\phi(\xi-s)}{s^{1+\alpha}}\,ds.
    \eeaa
    Differentiating under the integral sign and integrating by parts, we then obtain
    \beaa
        \phi'(\xi) \dalpha \phi(\xi) &=& \dfrac{\alpha}{2\Gamma(1-\alpha)} \int^{+\infty}_{0} \dfrac{1}{s^{1+\alpha}} \dfrac{d}{d\xi} [(\phi(\xi)-\phi(\xi-s))^2]\,ds \\
        &&\qquad + \dfrac{\alpha}{2\Gamma(1-\alpha)} \int^{+\infty}_{0} \dfrac{1}{s^{1+\alpha}} \dfrac{d}{ds} [(\phi(\xi)-\phi(\xi-s))^2]\,ds \\
        &=& \dfrac{\alpha}{2\Gamma(1-\alpha)} \dfrac{d}{d\xi} \int^{+\infty}_{0} \dfrac{(\phi(\xi)-\phi(\xi-s))^2}{s^{1+\alpha}}\,ds \\
        &&\qquad + \dfrac{\alpha(1+\alpha)}{2\Gamma(1-\alpha)} \int^{+\infty}_{0} \dfrac{(\phi(\xi)-\phi(\xi-s))^2}{s^{2+\alpha}}\,ds.
    \eeaa
    We integrate both sides over $\R$.
    By the dominated convergence theorem, the first term on the right-hand side contributes zero.
    We therefore obtain
    \be\label{ieq:dalpha-pos}
        \int_{\R} \phi'(\xi) \dalpha \phi(\xi)\,d\xi = \dfrac{\alpha(1+\alpha)}{2\Gamma(1-\alpha)} \int_{\R} \int^{+\infty}_{0} \dfrac{(\phi(\xi)-\phi(\xi-s))^2}{s^{2+\alpha}}\,ds\,d\xi >0.
    \ee

    Integrating both sides of \eqref{eq:pre-energy} over $\R$ and using \eqref{ieq:dalpha-pos} and \eqref{ass:nonlinearity}, we obtain
    \beaa
        \theta \dfrac{\alpha(1+\alpha)}{2\Gamma(1-\alpha)} \int_{\R} \int^{+\infty}_{0} \dfrac{(\phi(\xi)-\phi(\xi-s))^2}{s^{2+\alpha}}\,ds\,d\xi = F(1)> 0. 
    \eeaa
    Hence $\theta>0$.
\end{proof}
\begin{remark}\label{rem:speed}
    A similar argument shows that no asymptotic traveling wave with $c<0$ exists.
    Suppose that there exists $(c,\phi)\in(-\infty,0)\times C_b^2(\R)$ satisfying \eqref{ass:boundary_condition} and
    \beaa
        -\frac{|c|^\alpha}{\Gamma(1-\alpha)} 
        \int_{0}^{+\infty}
        \frac{\phi'(\xi+s)}{s^\alpha}\,ds
        =
        \phi''(\xi)+f(\phi(\xi))
        \qquad (\xi\in\R).
    \eeaa
    We multiply the equation by $\phi'$ and integrate over $\R$.
    A similar calculation then yields
    \beaa
        0\ge-\dfrac{|c|^\alpha\alpha(1+\alpha)}{2\Gamma(1-\alpha)} \int_{\R} \int^{+\infty}_{0} \dfrac{(\phi(\xi)-\phi(\xi+s))^2}{s^{2+\alpha}}\,ds\,d\xi = F(1)> 0,
    \eeaa
    which is a contradiction.
\end{remark}

\medskip
Next, we show that the profile lies strictly between the two stable equilibria.
For this, we use $\theta>0$ and Lemma~\ref{lem:max}.
\begin{lemma}\label{lem:apriori}
    If $(\theta,\phi)\in(0,+\infty)\times C_b^2(\R)$ is a solution of Equation~\eqref{eq:atw-theta} satisfying \eqref{ass:boundary_condition}, then $0<\phi<1$.
\end{lemma}
\begin{proof}
    Let $d_0(\xi) = \int_0^1 f'(t\phi(\xi))\,dt$ and $d_1(\xi)=\int_0^1f'(1+t(\phi(\xi)-1))\,dt$.
    Since $\phi$ is bounded, we have $d_j \in C_b(\R)$.
    Since $\phi$ solves Equation~\eqref{eq:atw-theta}, the functions $w_0=\phi$ and $w_1=1-\phi$ satisfy
    \beaa
        w_j'' -\theta \dalpha w_j +d_jw_j =0\qquad (j=0,1).
    \eeaa
    Set $H_j := \{\xi:w_j \ge 0\}$.
    If $\xi \not \in H_0$, then $\phi <0$, and hence $d_0 = f(\phi)/ \phi \le 0$.
    If $\xi \not \in H_1$, then $\phi >1$, and hence $d_1 = f(\phi)/ (\phi-1) \le 0$.
    Furthermore, \eqref{ass:boundary_condition} gives $\liminf_{\xi \to \pm \infty} w_j(\xi) \ge 0$.
    Thus Lemma~\ref{lem:max} implies that $w_0,w_1 >0$.
    Therefore, $0<\phi<1$.
\end{proof}

\medskip
We next establish the regularity of the profile and the decay of its derivatives at infinity.
These properties allow us to differentiate the profile equation and apply Lemmas~\ref{lem:min-sgn} and~\ref{lem:max} to the resulting equation.
\begin{lemma}\label{lem:regularity}
    If $(\theta,\phi)\in(0,+\infty)\times C_b^2(\R)$ is a solution of Equation~\eqref{eq:atw-theta}, then $\phi\in C_b^3(\R)$.
    Moreover, if $\phi' \ge 0$ and \eqref{ass:boundary_condition} holds, then $\phi'\in C_0^2(\R)$.
\end{lemma}
\begin{proof}
    By Lemma~\ref{lem:D-bdd}, $\dalpha\phi\in C_b^1(\R)$.
    Since $f\in C^1(\R)$ and $\phi,\phi'$ are bounded, we also have $f(\phi)\in C_b^1(\R)$.
    Thus $\phi''\in C_b^1(\R)$, and hence $\phi\in C_b^3(\R)$.

    Next, assume that $\phi'\ge0$ and \eqref{ass:boundary_condition} holds.
    Then $\phi'$ is uniformly continuous and integrable, so $\phi'\in C_0(\R)$.
    Moreover, Lemma~\ref{lem:Dalpha-phi-C0} gives $\dalpha\phi\in C_0(\R)$.
    Therefore, Equation~\eqref{eq:atw-theta} and the continuity of $f$ imply that
    \beaa
        \lim_{\xi\to-\infty}\phi''(\xi)=-f(0)=0,
        \qquad
        \lim_{\xi\to+\infty}\phi''(\xi)=-f(1)=0.
    \eeaa
    Thus $\phi''\in C_0(\R)$.
    Differentiating Equation~\eqref{eq:atw-theta} and using Lemma~\ref{lem:D-bdd}, we obtain
    \beaa
        \phi'''=\theta\dalpha(\phi')-f'(\phi)\phi'.
    \eeaa
    Since $\phi'\in C_0^1(\R)$, we also have $\dalpha(\phi')\in C_0(\R)$.
    Moreover, the boundedness of $f'(\phi)$ implies that $f'(\phi)\phi'\in C_0(\R)$.
    Hence $\phi'''\in C_0(\R)$, and therefore $\phi'\in C_0^2(\R)$.
\end{proof}

\medskip
We next establish the monotonicity of asymptotic traveling wave profiles.
\begin{proposition}\label{prop:mono}
    If $(\theta,\phi)\in(0,+\infty)\times C_b^2(\R)$ is a solution of Equation~\eqref{eq:atw-theta} satisfying \eqref{ass:boundary_condition}, then $\phi'>0$.
\end{proposition}
\begin{proof}
    For each $h>0$, let
    \beaa
        w_h(\xi):=\phi(\xi+h)-\phi(\xi),\qquad
        d_h(\xi):=\int^1_0 f'\bigl(\phi(\xi)+t w_h(\xi)\bigr)\,dt.
    \eeaa
    By Equation~\eqref{eq:atw-theta} and \eqref{ass:boundary_condition},
    we have
    \beaa
        w_h''-\theta \dalpha w_h+d_hw_h=0,
        \qquad w_h(\pm\infty)=0.
    \eeaa
    Moreover, $w_h\in C_b^2(\R)$ and $d_h\in C_b(\R)$.
    We first show that $w_h>0$ for every $h>0$.
    
    By \eqref{ass:nonlinearity}, there exists $\delta>0$ such that
    \beaa
        f'(u)<0
        \qquad
        \bigl(u\in[0,\delta]\cup[1-\delta,1]\bigr).
    \eeaa
    Moreover, by \eqref{ass:boundary_condition} and
    Lemma~\ref{lem:apriori}, we can choose $R_0>0$ sufficiently large
    so that
    \beaa
        1-\delta<\phi(\xi)<1
        \quad(\xi>R_0),
        \qquad
        0<\phi(\xi)<\delta
        \quad(\xi<-R_0).
    \eeaa
    For each $h>0$, set
    \beaa
        m(h):=\min_{\xi\in[-R_0,R_0]}w_h(\xi).
    \eeaa

    We first show that, for each $h>0$, $m(h)\ge0$ implies $w_h>0$. 
    Assume that $m(h)\ge0$, and set $H_h:=\{\xi\in\R\mid w_h(\xi)\ge0\}$.
    Then $[-R_0,R_0]\subset H_h$.
    Let $\xi\notin H_h$.
    If $\xi<-R_0$, then $0<\phi(\xi+h)<\phi(\xi)<\delta$.
    If $\xi>R_0$, then $\phi(\xi),\phi(\xi+h)\in(1-\delta,1)$.
    Thus, $d_h(\xi)<0$ in either case.
    By Lemma~\ref{lem:max}, we obtain either
    $w_h>0$ or $w_h\equiv0$.
    The latter would imply that $\phi$ is $h$-periodic,
    which contradicts \eqref{ass:boundary_condition}.
    Therefore, $w_h>0$, and in particular, $m(h)>0$.

    Next, we show that $m(h)>0$ for every $h>0$.
    Set
    \beaa
        \eta=\dfrac{1}{2} \left( 1- \max_{\xi\in[-R_0,R_0]} \phi(\xi) \right).
    \eeaa
    By Lemma~\ref{lem:apriori}, we have $\eta>0$.
    By \eqref{ass:boundary_condition}, there exists $R_1>0$ such that
    \beaa
        \phi(s)>1-\eta\qquad(s\ge R_1).
    \eeaa
    Thus, for $h\ge R_0+R_1$ and $\xi\in [-R_0,R_0]$,
    we have $\xi+h\ge R_1$, and hence
    \beaa
        w_h(\xi) =\phi(\xi+h)-\phi(\xi) \ge 1-\eta- \max_{\xi\in[-R_0,R_0]} \phi(\xi) =\eta.
    \eeaa
    Therefore, $m(h)\ge\eta$ for all $h\ge R_0+R_1$.

    Since $\phi\in C_b^1(\R)$, $m$ is Lipschitz continuous.
    Suppose that $m(h)\le 0$ for some $h>0$.
    Then, by continuity and the preceding estimate,
    there exists $h^*>0$ such that $m(h^*)=0$.
    However, the previous argument implies that
    $m(h^*)>0$, a contradiction.
    Therefore, $m(h)>0$ for every $h>0$.
    
    Thus, $w_h>0$ for every $h>0$.
    Finally, $\phi'=\lim_{h\downarrow0}w_h/h\ge0$.
    By Lemma~\ref{lem:regularity}, we have $\phi'\in C_b^2(\R)\cap C_0^1(\R)$.
    Differentiating Equation~\eqref{eq:atw-theta}, we obtain
    \beaa
        (\phi')''-\theta\dalpha \phi'+f'(\phi)\phi'=0.
    \eeaa
    Since $\phi'\ge0$, Lemma~\ref{lem:max} implies that
    either $\phi'>0$ or $\phi'\equiv0$.
    The latter contradicts \eqref{ass:boundary_condition}.
    Therefore, $\phi'>0$ on $\R$.
\end{proof}

\medskip
Finally, we prove that the wave speed is unique and that the profile is unique up to translation.
\begin{proposition}\label{prop:uni}
    Let $(\theta_j,\phi_j)\in(0,+\infty)\times C_b^2(\R)$, $j=1,2$, be solutions of Equation~\eqref{eq:atw-theta} satisfying \eqref{ass:boundary_condition}.
    Then $\theta_1=\theta_2$.
    Moreover, there exists $h\in\R$ such that $\phi_1(\xi)=\phi_2(\xi+h)$ for all $\xi\in\R$.
\end{proposition}
\begin{proof}
    Without loss of generality, assume that $0<\theta_2\le\theta_1$.
    For $h\in\R$, define
    \beaa
        w_h(\xi):=\phi_2(\xi+h)-\phi_1(\xi),\qquad
        d_h(\xi):=\int_0^1 f'\bigl(\phi_1(\xi)+t w_h(\xi)\bigr)\,dt.
    \eeaa
    Proposition~\ref{prop:mono} implies that $\dalpha\phi_1>0$.
    Subtracting the two profile equations, we obtain
    \be\label{eq:uni-diff}
        w_h''-\theta_2\dalpha w_h+d_hw_h
        =-(\theta_1-\theta_2)\dalpha\phi_1\le0.
    \ee
    Moreover, $d_h\in C_b(\R)$ and $w_h(\pm\infty)=0$.

    Since $\phi_2(-\infty)=0<\phi_1(0)$, we can fix $h^-\in\R$ such that
    \beaa
        \phi_2(h^-)<\phi_1(0).
    \eeaa
    Choose $\delta\in(0,1/2)$ such that
    \beaa
        f'(u)<0
        \qquad
        \bigl(u\in[0,\delta]\cup[1-\delta,1]\bigr).
    \eeaa
    By \eqref{ass:boundary_condition}, we can choose $R_0>0$ sufficiently large so that
    \beaa
        \phi_1(\xi)\in(0,\delta)
        \quad(\xi<-R_0),\qquad
        \phi_1(\xi),\ \phi_2(\xi+h)\in(1-\delta,1)
        \quad(\xi>R_0,\ h\ge h^-).
    \eeaa

    Set
    \beaa
        m(h):=\min_{\xi\in[-R_0,R_0]}w_h(\xi)
        \qquad(h\ge h^-).
    \eeaa
    Since $0\in[-R_0,R_0]$, we have
    \beaa
        m(h^-)\le\phi_2(h^-)-\phi_1(0)<0.
    \eeaa
    On the other hand, set
    \beaa
        \eta=\dfrac{1}{2}
        \left(1-\max_{\xi\in[-R_0,R_0]}\phi_1(\xi)\right)>0.
    \eeaa
    By \eqref{ass:boundary_condition}, there exists $R_1>0$ such that
    $\phi_2(s)>1-\eta$ for all $s\ge R_1$.
    Choose $h^+>\max\{h^-,R_0+R_1\}$.
    Then, for every $\xi\in[-R_0,R_0]$, we have
    \beaa
        w_{h^+}(\xi)
        \ge 1-\eta-\max_{\xi\in[-R_0,R_0]}\phi_1(\xi)
        =\eta.
    \eeaa
    Thus, $m(h^+)>0$.

    Since $\phi_2\in C_b^1(\R)$, $m$ is Lipschitz continuous.
    The intermediate value theorem therefore gives $h^*\in[h^-,h^+]$ such that $m(h^*)=0$.
    Hence $w_{h^*}\ge0$ on $[-R_0,R_0]$.
    Moreover, $w_{h^*}(\xi^*)=0$ for some $\xi^*\in[-R_0,R_0]$.

    Set $H:=\{\xi\in\R\mid w_{h^*}(\xi)\ge0\}$.
    Then $[-R_0,R_0]\subset H$.
    Let $\xi\notin H$.
    If $\xi<-R_0$, then
    \beaa
        0<\phi_2(\xi+h^*)<\phi_1(\xi)<\delta.
    \eeaa
    If $\xi>R_0$, then both profile values lie in $(1-\delta,1)$.
    Thus, $d_{h^*}(\xi)<0$ in either case.
    Since $w_{h^*}(\xi^*)=0$, Lemma~\ref{lem:max} implies that $w_{h^*}\equiv0$.

    Substituting this into \eqref{eq:uni-diff}, we obtain
    \beaa
        (\theta_1-\theta_2)\dalpha\phi_1=0.
    \eeaa
    Since $\dalpha\phi_1>0$, we conclude that $\theta_1=\theta_2$.
    Moreover, $w_{h^*}\equiv0$ gives $\phi_1(\xi)=\phi_2(\xi+h^*)$.
\end{proof}

\medskip
We have shown that asymptotic traveling waves satisfying \eqref{ass:boundary_condition} have positive speed and strictly increasing profiles.
The speed is unique, and the profile is unique up to translation.
In the next section, we prove the existence of such waves.

\bigskip
\section{Existence of asymptotic traveling waves}\label{sec:proof}

In this section, we prove the existence of asymptotic traveling waves.
As outlined above, we use a continuation method.
We start from the stationary front in the balanced case described in Section~\ref{subsec:stasol}.
We continuously deform the nonlinearity into $f$ and construct a continuous family of solutions starting from this front.

First, given $f$ satisfying \eqref{ass:nonlinearity}, we construct $g\in C^1(\R)$ such that
\beaa
    \int_0^1g(u)\,du=0.
\eeaa
Choose $\chi\in C_c^\infty(\R)$ satisfying $\chi\ge0$, $\chi\not\equiv0$, and $\operatorname{supp}\chi\subset[0,a]$.
Define
\beaa
    g(u)
    =
    f(u)+\lambda\chi(u)f(u),
    \qquad
    \lambda
    =
    -\dfrac{\int_0^1 f(u)\,du}{\int_0^a\chi(u)f(u)\,du}.
\eeaa
By \eqref{ass:nonlinearity}, we have $\lambda>0$.
The function $g$ has the same zeros and sign as $f$.
We then define
\beaa
    f_\eps(u)
    :=
    (1-\eps)g(u)+\eps f(u),
    \qquad
    \eps\in\R.
\eeaa
Since $f_\eps$ satisfies \eqref{ass:nonlinearity} for each $\eps\in(0,1]$, the arguments in Section~\ref{sec:mono-uni} also apply to $f_\eps$.

In what follows, we restrict $\eps$ to $[0,1]$.
Fix $\tilde a\in(a,1)$ and consider the normalized problem
\be\label{eq:problem-atw}
    \begin{cases}
        \phi''-\theta\dalpha\phi+f_\eps(\phi)=0,
        \\[1mm]
        \phi(-\infty)=0,
        \qquad
        \phi(+\infty)=1,
        \qquad
        \phi(0)=\tilde a.
    \end{cases}
\ee
The condition $\phi(0)=\tilde a$ fixes the translation.
For $\eps=0$, Proposition~\ref{prop:stationary-front} gives a monotone stationary front $\phi_0$.
By a suitable translation, we may assume that $\phi_0(0)=\tilde a$.
Then $(\theta_0,\phi_0)=(0,\phi_0)$ solves Problem~\eqref{eq:problem-atw} at $\eps=0$.

We now prove the existence of solutions to Problem~\eqref{eq:problem-atw} for $\eps\in(0,1]$.
In Section~\ref{subsec:limsol}, we show that solutions exist for every $\eps\in[0,\eps_0)$ for some $\eps_0>0$.
By passing to the limit, we also obtain a solution at $\eps=\eps_0$.
In Section~\ref{subsec:extension}, we assume that a solution exists at $\eps=\eps_0\in(0,1)$.
We show that this solution extends to a continuous family near $\eps_0$.
Using these results, we prove Theorem~\ref{thm:main} in Section~\ref{subsec:proof}.

\medskip
\subsection{Local existence and passage to the limit}\label{subsec:limsol}

We first show that Problem~\eqref{eq:problem-atw} has a continuous family of solutions for sufficiently small $\eps>0$.
\begin{proposition}\label{prop:local-zero}
    There exists $\eps_0>0$ such that, for every $\eps\in[0,\eps_0)$, Problem \eqref{eq:problem-atw} has a solution $(\theta_\eps,\phi_\eps)\in\R\times C_b^2(\R)$.
    Moreover, the map $\eps\mapsto(\theta_\eps,\phi_\eps)$ is of class $C^1$.
\end{proposition}
\begin{proof}
    For $\sigma\in\R$ and $z=(\psi,\tau)\in C_0^2(\R)\times\R$, define
    \beaa
        H_0(\sigma,z)
        :=
        \left(
        (\phi_0+\psi)''-\tau\dalpha(\phi_0+\psi)+f_\sigma(\phi_0+\psi),
        \psi(0)
        \right).
    \eeaa
    Then $H_0:\R\times(C_0^2(\R)\times\R)\to C_0(\R)\times\R$ is of class $C^1$.
    Moreover, $H_0(0,0)=(0,0)$, and
    \beaa
        DH_0
        :=
        \frac{\partial H_0}{\partial z}(0,0)
        =
        \begin{bmatrix}
            L_0 & -\dalpha\phi_0\\
            E_0 & 0
        \end{bmatrix},
        \qquad
        E_0\psi:=\psi(0).
    \eeaa

    We show that $DH_0:C_0^2(\R)\times\R\to C_0(\R)\times\R$ is an isomorphism.
    The conclusion will then follow from the implicit function theorem.
    For any $(A,b)\in C_0(\R)\times\R$, consider
    \beaa
        L_0\psi=A+\tau\dalpha\phi_0,
        \qquad
        \psi(0)=b.
    \eeaa
    By Lemma~\ref{lem:Dalpha-phi-C0}, we have $\dalpha\phi_0\in C_0(\R)$.
    Thus, Lemma~\ref{lem:solvability-L0} gives the following necessary and sufficient condition for solvability in $C_0^2(\R)$:
    \be\label{solve:tau}
        \int_{\R}(A+\tau\dalpha\phi_0)\phi_0'\,d\xi=0.
    \ee
    The solution is unique whenever this condition holds.
    Since $\phi_0'>0$, we have $\dalpha\phi_0>0$.
    Hence \eqref{solve:tau} uniquely determines $\tau$.
    For this $\tau$, the solution $\psi\in C_0^2(\R)$ is also uniquely determined.
    Therefore, $DH_0$ is a bounded linear bijection.
    By the bounded inverse theorem, it is an isomorphism.

    The implicit function theorem now yields $\eps_0>0$ and a $C^1$ map
    \beaa
        (-\eps_0,\eps_0)\ni\eps
        \longmapsto
        (\psi_\eps,\tau_\eps)\in C_0^2(\R)\times\R
    \eeaa
    satisfying $H_0(\eps,(\psi_\eps,\tau_\eps))=(0,0)$ and $(\psi_0,\tau_0)=(0,0)$.
    For $\eps\in[0,\eps_0)$, set
    \beaa
        \phi_\eps:=\phi_0+\psi_\eps,
        \qquad
        \theta_\eps:=\tau_\eps.
    \eeaa
    The profile equation and normalization in Problem~\eqref{eq:problem-atw} follow from $H_0(\eps,(\psi_\eps,\tau_\eps))=(0,0)$.
    The conditions \eqref{ass:boundary_condition} follow from $\psi_\eps\in C_0^2(\R)$.
    This proves the proposition.
\end{proof}

\begin{remark}\label{rem:speed-zero}
    For use in Section~\ref{sec:num}, we derive the first-order expansion of $\theta_\eps$ at $\eps=0$.
    A dot denotes differentiation with respect to $\eps$.
    We differentiate the first equation and the normalization condition in Problem~\eqref{eq:problem-atw} with respect to $\eps$ at $\eps=0$.
    Since $\theta_0=0$, we obtain
    \beaa
        L_0\dot\psi_0
        =
        \dot\theta_0\dalpha\phi_0-(f-g)(\phi_0),
        \qquad
        \dot\psi_0(0)=0.
    \eeaa
    Since $\dot\psi_0\in C_0^2(\R)$, the solvability condition in Lemma~\ref{lem:solvability-L0} yields
    \beaa
        \dot\theta_0
        =
        \dint_{\R}(f-g)(\phi_0)\phi_0'\,d\xi
        \Big/
        \dint_{\R}\dalpha\phi_0\,\phi_0'\,d\xi
        =
        \dint_0^1 f(u)\,du
        \Big/
        \dint_{\R}\dalpha\phi_0\,\phi_0'\,d\xi.
    \eeaa
    Hence $\theta_\eps=\dot\theta_0\eps+o(\eps)$ as $\eps\to +0$.
    Since $\theta_\eps=c^\alpha$, we obtain
    \beaa
        c=(\dot\theta_0)^{1/\alpha}\eps^{1/\alpha}
        +o(\eps^{1/\alpha})
        \qquad(\eps\to +0).
    \eeaa
\end{remark}

\medskip
We now construct a solution at $\eps=\eps_0$.
In the following argument and in Section~\ref{subsec:extension}, fix $\eps_0\in(0,1]$.
Assume that, for every $\eps\in[0,\eps_0)$, Problem~\eqref{eq:problem-atw} has a solution $(\theta_\eps,\phi_\eps)\in\R\times C_b^2(\R)$.
For $\eps>0$, Propositions~\ref{prop:posi-speed} and~\ref{prop:mono} imply that $\phi_\eps'>0$.
The same holds at $\eps=0$ by Proposition~\ref{prop:stationary-front}.

To construct a solution at $\eps=\eps_0$, take any sequence $\{\eps_n\}\subset[0,\eps_0)$ with $\eps_n\nearrow\eps_0$.
Write $(\theta_n,\phi_n)$ for the corresponding solutions of Problem~\eqref{eq:problem-atw}.

We first establish a uniform bound on $\theta_n$.
This bound will be used when passing to the limit.
\begin{lemma}\label{lem:theta-boundness}
    The sequence $\{\theta_n\}$ is bounded.
\end{lemma}
\begin{proof}
    Suppose that $\{\theta_n\}$ is unbounded.
    Since $\eps_n>0$ for all sufficiently large $n$, Proposition~\ref{prop:posi-speed} gives $\theta_n>0$ for such $n$.
    Passing to a subsequence, we may assume that $\theta_n>0$ and $\theta_n\to+\infty$.
    Set $p_n:=\phi_n'$.
    By Lemma~\ref{lem:regularity}, $p_n\in C_0^2(\R)$.
    Since $p_n>0$, it attains its global maximum at some $\xi_n\in\R$.

    Set $M_n:=p_n(\xi_n)$.
    We differentiate the first equation in Problem~\eqref{eq:problem-atw} using Lemma~\ref{lem:D-bdd}.
    At $\xi=\xi_n$, this gives
    \beaa
        \theta_n\dalpha p_n(\xi_n)
        =
        p_n''(\xi_n)+f_{\eps_n}'(\phi_n(\xi_n))M_n
        \le LM_n.
    \eeaa
    Here, $L>0$ is a constant such that $|f_\eps'(u)|\le L$ on $[0,\eps_0]\times[0,1]$.
    Set $C:=L\Gamma(1-\alpha)/\alpha$.
    Applying Lemma~\ref{lem:min-sgn} to $-p_n$ and using $\|p_n\|_{L^1(\R)}=1$, we obtain
    \beaa
        \frac{C}{\theta_n}M_n
        \ge
        \frac{1}{R^{1+\alpha}}
        \int_0^R\bigl(M_n-p_n(\xi_n-s)\bigr)\,ds
        \ge
        \frac{RM_n-1}{R^{1+\alpha}}
    \eeaa
    for any $R>0$.
    Taking $R=2/M_n$, we obtain $M_n^\alpha\le 2^{1+\alpha}C/\theta_n\to0$.
    Thus, $M_n\to0$.

    On the other hand, fix $0<u_-<u_+<a$ and set $U:=[u_-,u_+]$.
    By the monotonicity of $\phi_n$ and \eqref{ass:boundary_condition}, there exist $a_n<b_n$ such that
    $\phi_n(a_n)=u_-$ and $\phi_n(b_n)=u_+$.
    Hence
    \beaa
        u_+-u_-
        =\int_{a_n}^{b_n}p_n(\xi)\,d\xi
        \le M_n(b_n-a_n).
    \eeaa
    The mean value theorem gives $\tilde\xi_n\in(a_n,b_n)$ such that
    \beaa
        p_n(b_n)-p_n(a_n)
        =
        \phi_n''(\tilde\xi_n)(b_n-a_n).
    \eeaa
    Since $0<p_n\le M_n$, we have
    \beaa
        |\phi_n''(\tilde\xi_n)|
        \le\frac{M_n}{b_n-a_n}
        \le\frac{M_n^2}{u_+-u_-}
        \longrightarrow0.
    \eeaa
    Moreover, by continuity and \eqref{ass:nonlinearity}, there exists $\kappa>0$ such that
    \beaa
        f_\eps(u)\le-\kappa
        \qquad
        \bigl(\eps\in[0,\eps_0],\ u\in U\bigr).
    \eeaa
    Since $\phi_n(\tilde\xi_n)\in U$, the first equation in Problem~\eqref{eq:problem-atw} gives, for all sufficiently large $n$,
    \beaa
        0<
        \theta_n(\dalpha\phi_n)(\tilde\xi_n)
        =
        \phi_n''(\tilde\xi_n)
        +f_{\eps_n}(\phi_n(\tilde\xi_n))
        \le
        \frac{M_n^2}{u_+-u_-}-\kappa
        <0,
    \eeaa
    which is a contradiction.
    Therefore, $\{\theta_n\}$ is bounded.
\end{proof}

We next prove the following lemma to extract a convergent subsequence from $\{\phi_n\}$.
\begin{lemma}\label{lem:uni-bdd}
    The sequence $\{\phi_n\}$ is uniformly bounded in $C_b^3(\R)$.
\end{lemma}
\begin{proof}
    For $\zeta_1,\zeta_2\ge0$, we write $\zeta_1\lesssim\zeta_2$ if $\zeta_1\le C\zeta_2$ for some constant $C$ independent of $n$.

    The functions $f_\eps$ and $f_\eps'$ are bounded on $[0,\eps_0]\times[0,1]$.
    Together with Lemmas~\ref{lem:apriori} and~\ref{lem:theta-boundness}, this gives
    \beaa
        |\theta_n|
        +\|f_{\eps_n}(\phi_n)\|_{L^\infty}
        +\|f_{\eps_n}'(\phi_n)\|_{L^\infty}
        \lesssim1.
    \eeaa
    The profile equation can be written as
    \be\label{eq:phi-n}
        \phi_n''
        =\theta_n\dalpha\phi_n-f_{\eps_n}(\phi_n).
    \ee
    By \eqref{ineq:Dalpha-est} and Landau's inequality~\cite{Hille1972}, we obtain
    \beaa
        \|\phi_n''\|_{L^\infty}
        \lesssim1+\|\phi_n'\|_{L^\infty}
        \lesssim1+\sqrt{\|\phi_n''\|_{L^\infty}}.
    \eeaa
    Thus, $\|\phi_n''\|_{L^\infty}\lesssim1$.
    Applying Landau's inequality again gives $\|\phi_n'\|_{L^\infty}\lesssim1$.
    Hence $\|\phi_n\|_{C_b^2}\lesssim1$.
    Finally, differentiating \eqref{eq:phi-n} using Lemma~\ref{lem:D-bdd}, we obtain
    \beaa
        \phi_n'''
        =\theta_n\dalpha(\phi_n')
        -f_{\eps_n}'(\phi_n)\phi_n'.
    \eeaa
    The estimate in the same lemma gives
    \beaa
        \|\phi_n'''\|_{L^\infty}
        \lesssim\|\phi_n\|_{C_b^2}
        \lesssim1.
    \eeaa
    Therefore, $\{\phi_n\}$ is uniformly bounded in $C_b^3(\R)$.
\end{proof}

\medskip
Lemma~\ref{lem:uni-bdd} gives relative compactness of $\{\phi_n\}$ in $C_{\mathrm{loc}}^2(\R)$.
We obtain the profile equation by passing to the limit along a subsequence.
We then verify the conditions at infinity to construct a solution at $\eps=\eps_0$.
\begin{proposition}\label{prop:limit-solution}
    There exists $(\theta_{\eps_0},\phi_{\eps_0})\in(0,+\infty)\times C_b^2(\R)$ which solves Problem~\eqref{eq:problem-atw} at $\eps=\eps_0$.
\end{proposition}
\begin{proof}
    By Lemmas~\ref{lem:theta-boundness} and~\ref{lem:uni-bdd}, the sequences $\{\theta_n\}$ and $\{\phi_n\}$ are bounded in $\R$ and $C_b^3(\R)$, respectively.
    By the Ascoli--Arzel\`a theorem and a diagonal argument, we may pass to a subsequence such that
    \beaa
        \theta_n\to\bar\theta\ge0,\qquad
        \phi_n\to\ovphi\quad\text{in }C_{\mathrm{loc}}^2(\R).
    \eeaa
    Then $\ovphi\in C_b^2(\R)$ and
    \beaa
        0\le\ovphi\le1,\qquad
        \ovphi(0)=\tilde a,\qquad
        \ovphi'\ge0.
    \eeaa
    Using Lemma~\ref{lem:Dalpha-local-conv}, we pass to the limit in the first equation of Problem~\eqref{eq:problem-atw} and obtain
    \be\label{eq:limit-equation}
        \bar\theta\dalpha\ovphi-\ovphi''
        =f_{\eps_0}(\ovphi)
        \qquad\text{in }\R.
    \ee
    Finally, we verify \eqref{ass:boundary_condition} for $\ovphi$.

    Since $\ovphi$ is bounded and nondecreasing, the limits
    $L^\pm:=\lim_{\xi\to\pm\infty}\ovphi(\xi)$ exist.
    As in the proof of Lemma~\ref{lem:Dalpha-phi-C0}, we have $\ovphi'\in C_0(\R)$.
    The lemma also gives $\dalpha\ovphi\in C_0(\R)$.
    Integrating \eqref{eq:limit-equation} over $[\xi,\xi+1]$ and letting $\xi\to\pm\infty$, we obtain $f_{\eps_0}(L^\pm)=0$.
    Since $L^-\le\tilde a\le L^+$ and $a<\tilde a<1$, we have $L^+=1$ and $L^-\in\{0,a\}$.

    For both $v=\phi_n$ and $v=\ovphi$, we have $v'\in L^1(\R)\cap C_0(\R)$.
    Lemma~\ref{lem:Dalpha-phi-C0} also gives $\dalpha v\in C_0(\R)$.
    We multiply the corresponding equation by $v'$ and integrate over $\R$.
    Since $\int_\R v''v'\,d\xi=0$, this gives
    \be\label{eq:energy}
        \theta\int_\R(\dalpha v)v'\,d\xi
        =
        \int_{v(-\infty)}^{v(+\infty)}f_\eps(u)\,du,
    \ee
    where $(\theta,\eps)$ are the parameters corresponding to $v$.
    By monotonicity, $(\dalpha\phi_n)\phi_n'\ge0$.
    Fatou's lemma and \eqref{eq:energy} therefore yield
    \beaa
        \bar\theta\int_\R(\dalpha\ovphi)\ovphi'\,d\xi
        \le
        \liminf_{n\to\infty}
        \theta_n\int_\R(\dalpha\phi_n)\phi_n'\,d\xi
        =
        \int_0^1f_{\eps_0}(u)\,du.
    \eeaa
    Suppose that $L^-=a$.
    By \eqref{eq:energy}, the left-hand side of the preceding inequality equals $\int_a^1f_{\eps_0}(u)\,du$.
    Hence $0\le\int_0^a f_{\eps_0}(u)\,du$, contradicting $f_{\eps_0}<0$ on $(0,a)$.
    Therefore, $L^-=0$.

    Thus, $(\theta_{\eps_0},\phi_{\eps_0}):=(\bar\theta,\ovphi)$ solves Problem~\eqref{eq:problem-atw}.
    Propositions~\ref{prop:posi-speed} and~\ref{prop:mono} imply that $\theta_{\eps_0}>0$ and $\phi_{\eps_0}'>0$.
\end{proof}

\medskip
We have thus extended the interval of existence from $(0,\eps_0)$ to $(0,\eps_0]$.
If $\eps_0<1$, local existence near $\eps=\eps_0$ would extend this interval further.
We therefore seek a result analogous to Proposition~\ref{prop:local-zero}.
We establish the necessary lemmas in the next subsection.

\medskip
\subsection{Continuation of solutions}\label{subsec:extension}
In this subsection,  we assume that $\eps_0<1$.
We prove the following proposition.
\begin{proposition}\label{prop:extension}
    There exists $\eta>0$ such that, for every $\eps\in[\eps_0,\eps_0+\eta)$,
    Problem \eqref{eq:problem-atw} has a solution $(\theta_\eps,\phi_\eps)\in(0,+\infty)\times C_b^2(\R)$.
    Moreover, the map $\eps\mapsto(\theta_\eps,\phi_\eps)$ is of class $C^1$.
\end{proposition}

\medskip
The proof uses the implicit function theorem, as in Proposition~\ref{prop:local-zero}.
To this end, we first establish several lemmas concerning the linearized operator
\beaa
    L_{\eps_0}:C_0^2(\R)\to C_0(\R),
    \qquad
    L_{\eps_0}h
    :=
    h''-\theta_{\eps_0}\dalpha h
    +f_{\eps_0}'(\phi_{\eps_0})h.
\eeaa
\begin{lemma}\label{lem:fred}
    $L_{\eps_0}$ is Fredholm with index $0$.
\end{lemma}
\begin{proof}
    Let $\theta:=\theta_{\eps_0}>0$ and $b(\xi):=f_{\eps_0}'(\phi_{\eps_0}(\xi))$.
    By \eqref{ass:nonlinearity} and the continuity of $f_{\eps_0}'$, $b\in C_b(\R)$.
    Moreover, there exist $\mu,R>0$ such that $b(\xi)\le-\mu$ for $|\xi|\ge R$.
    Define $b_0(\xi):=\min\{b(\xi),-\mu\}$.
    Then $b_0\in C_b(\R)$, $b_0\le-\mu$, and $b-b_0$ is supported in $[-R,R]$.

    Choose $\kappa>0$ sufficiently large so that $\kappa^2\ge\mu$, and define
    \beaa
        B_s h&:=&h''-\theta\dalpha h
              +\bigl((1-s)(-\kappa^2)+sb_0\bigr)h
        \qquad (0\le s\le 1) \\
        K h&:=& (b-b_0)h.
    \eeaa
    By Lemma~\ref{lem:D-bdd}, $L_{\eps_0}$ and each $B_s$ are bounded linear operators from
    $C_0^2(\R)$ to $C_0(\R)$.
    Moreover,
    \be\label{Leps0}
        L_{\eps_0} = B_1 + K
    \ee
    holds.
    We prove the result by showing that $B_1$ is an isomorphism and $K$ is compact.

    We begin by showing that $B_0$ is surjective.
    We apply \cite[Proposition~3.5]{Ishii2025} with $c=\theta^{1/\alpha}$.
    For sufficiently large $\kappa>0$, there exists $G\in L^1(\R)$ such that, for any $A\in C_0(\R)$,
    \beaa
        h:=-G*A\in C_b^2(\R),
        \qquad
        B_0h=A.
    \eeaa
    Since $G\in L^1(\R)$ and $A\in C_0(\R)$, we have $h\in C_0(\R)$.
    Taylor's theorem gives, for any $r>0$,
    \beaa
        |h'(\xi)|
        \le
        \frac{|h(\xi+r)|+|h(\xi)|}{r}
        +\frac{r}{2}\|h''\|_{L^\infty}.
    \eeaa
    First letting $\xi\to\pm\infty$ and then $r\to +0$, we obtain $h'\in C_0(\R)$.
    Lemma~\ref{lem:D-bdd}~(iii) then gives $\dalpha h\in C_0(\R)$.
    The equation $B_0h=A$ then gives $h''\in C_0(\R)$.
    Thus, $B_0:C_0^2(\R)\to C_0(\R)$ is surjective.

    We next show that there exists a constant $C>0$, independent of $s$, such that
    \be\label{ineq:hom-uni-est}
        \|h\|_{C_b^2}\le C\|B_sh\|_{C_b}
    \ee
    for all $s\in[0,1]$ and $h\in C_0^2(\R)$.
    Set $A:=B_sh$.
    Suppose that $h\not\equiv0$.
    Replacing $(h,A)$ by $(-h,-A)$ if necessary, we may choose $\xi_*\in\R$ such that $h(\xi_*)=\|h\|_{L^\infty}>0$.
    At this point, $h''(\xi_*)\le0$.
    Lemma~\ref{lem:min-sgn} also gives $\dalpha h(\xi_*)\ge0$.
    Since $(1-s)(-\kappa^2)+sb_0\le-\mu$, we obtain
    \beaa
        A(\xi_*)\le-\mu h(\xi_*)=-\mu\|h\|_{L^\infty}.
    \eeaa
    Hence
    \be\label{ineq:hom-est}
        \|h\|_{L^\infty}\le\mu^{-1}\|A\|_{L^\infty}.
    \ee

    For $\zeta_1,\zeta_2\ge0$, we write $\zeta_1\lesssim\zeta_2$ if $\zeta_1\le C\zeta_2$ for some $C>0$ independent of $s$ and $h$.
    Using $B_sh=A$, Lemma~\ref{lem:D-bdd}, Landau's inequality, and \eqref{ineq:hom-est}, we obtain
    \beaa
        \|h''\|_{L^\infty}
        \lesssim
        \|A\|_{L^\infty}+\|h'\|_{L^\infty}
        \lesssim
        \|A\|_{L^\infty}
        +\|A\|_{L^\infty}^{1/2}\|h''\|_{L^\infty}^{1/2}.
    \eeaa
    Thus, $\|h''\|_{L^\infty}\lesssim\|A\|_{L^\infty}$.
    Applying Landau's inequality again gives $\|h'\|_{L^\infty}\lesssim\|A\|_{L^\infty}$.
    Together with \eqref{ineq:hom-est}, these bounds prove \eqref{ineq:hom-uni-est}.

    We have $B_s=(1-s)B_0+sB_1$.
    By \eqref{ineq:hom-uni-est}, the surjectivity of $B_0$, and \cite[Theorem~5.2]{GT}, $B_1$ is also surjective.
    Moreover, \eqref{ineq:hom-uni-est} implies that $B_1$ is injective and its inverse is bounded.
    Hence $B_1:C_0^2(\R)\to C_0(\R)$ is an isomorphism.

    Finally, we show that $K:C_0^2(\R)\to C_0(\R)$ is compact.
    Let $\{h_n\}$ be any bounded sequence in $C_0^2(\R)$.
    By the Ascoli--Arzel\`a theorem, we may pass to a subsequence that converges uniformly on $[-R,R]$.
    Since $b-b_0$ is supported in this interval, we have
    \beaa
        \|Kh_n-Kh_m\|_{L^\infty}
        \le
        \|b-b_0\|_{L^\infty}
        \|h_n-h_m\|_{C([-R,R])}
        \longrightarrow0
        \qquad(n,m\to\infty).
    \eeaa
    By the completeness of $C_0(\R)$, $\{Kh_n\}$ converges in $C_0(\R)$.
    Thus, $K$ is compact.

    By \eqref{Leps0}, $L_{\eps_0}$ is a compact perturbation of the isomorphism $B_1$.
    Compact perturbations preserve the Fredholm property and index.
    Hence, the desired assertion is obtained.
\end{proof}

\medskip
We next study the kernel and range of $L_{\eps_0}$.
The proof follows the idea in \cite[Lemma~3.9]{Bates2006}.
\begin{lemma}\label{lem:kernel-L}
    We have
    $\ker L_{\eps_0}=\operatorname{span}\{\phi_{\eps_0}'\}$.
    Moreover, if $A\in C_0(\R)\setminus\{0\}$ does not change sign, then $A\notin\operatorname{Ran}L_{\eps_0}$.
\end{lemma}
\begin{proof}
    Set $p:=\phi_{\eps_0}'$.
    Lemma~\ref{lem:regularity} and Propositions~\ref{prop:posi-speed} and~\ref{prop:mono} give $p>0$, $p\in C_0^2(\R)$, and $\theta_{\eps_0}>0$.
    Differentiating the first equation in Problem~\eqref{eq:problem-atw}, we obtain $L_{\eps_0}p=0$.
    Hence $\operatorname{span}\{p\}\subset\ker L_{\eps_0}$.

    To prove the reverse inclusion, let $q\in C_0^2(\R)$ satisfy $L_{\eps_0}q=0$.
    We show that $q$ is a constant multiple of $p$.

    Choose $R>0$ sufficiently large that
    \beaa
        f_{\eps_0}'(\phi_{\eps_0}(\xi))\le0
        \qquad(|\xi|>R).
    \eeaa
    The coefficient $f_{\eps_0}'(\phi_{\eps_0})$ is also bounded.
    Define
    \be\label{lambda-w}
        \lambda_0
        :=\max_{\xi\in[-R,R]}
        \left(-\frac{q(\xi)}{p(\xi)}\right),
        \qquad
        w:=q+\lambda_0p.
    \ee
    Since $p>0$, $\lambda_0$ is finite.
    By construction, $w\ge0$ on $[-R,R]$.
    Moreover, $w(\xi_0)=0$ for some $\xi_0\in[-R,R]$.
    We have $L_{\eps_0}w=0$ and $w(\pm\infty)=0$.
    Applying Lemma~\ref{lem:max} with $H=[-R,R]$, we obtain either $w>0$ or $w\equiv0$.
    Since $w(\xi_0)=0$, we have $w\equiv0$, and hence $q=-\lambda_0p$.
    Therefore, $\ker L_{\eps_0}=\operatorname{span}\{\phi_{\eps_0}'\}$.

    Finally, let $A\in C_0(\R)\setminus\{0\}$ satisfy either $A\le0$ or $A\ge0$.
    Suppose that $A\in\operatorname{Ran}L_{\eps_0}$.
    Then there exists $h\in C_0^2(\R)$ such that $L_{\eps_0}h=A$.
    If $A\le0$, define $\lambda_0$ and $w$ by replacing $q$ with $h$ in \eqref{lambda-w}.
    Then $L_{\eps_0}w=A\le0$, and $w\ge0$ on $[-R,R]$.
    Moreover, $w$ vanishes at some point in this interval.
    Applying Lemma~\ref{lem:max} again, we obtain $w\equiv0$.
    This contradicts $A\not\equiv0$.

    If $A\ge0$, applying the same argument to $L_{\eps_0}(-h)=-A$ gives a contradiction.
    Thus, $A\notin\operatorname{Ran}L_{\eps_0}$.
\end{proof}

\medskip
By Lemmas~\ref{lem:fred} and~\ref{lem:kernel-L}, $\operatorname{Ran}L_{\eps_0}$ is a closed subspace of $C_0(\R)$ of codimension $1$.
Hence there exists a nonzero bounded linear functional $\Lambda_{\eps_0}\in(C_0(\R))^*$ satisfying $\operatorname{Ran}L_{\eps_0}=\ker\Lambda_{\eps_0}$.

The functional $\Lambda_{\eps_0}$ has the following properties.
\begin{lemma}\label{lem:functional}
    Replacing $\Lambda_{\eps_0}$ by $-\Lambda_{\eps_0}$, if necessary, 
    we may assume that $\Lambda_{\eps_0}$ is a positive functional.
    Moreover, $\Lambda_{\eps_0}\bigl(\dalpha\phi_{\eps_0}\bigr)>0$.
\end{lemma}
\begin{proof}
    Set $\rho(\xi):=e^{-\xi^2}$.
    Since $-\rho\in C_0(\R)$ and $-\rho<0$, Lemma~\ref{lem:kernel-L} gives $-\rho\notin\operatorname{Ran}L_{\eps_0}$.
    Hence $\Lambda_{\eps_0}(\rho)\ne0$.
    Changing the sign of $\Lambda_{\eps_0}$ if necessary, we may assume that $\Lambda_{\eps_0}(\rho)>0$.

    Suppose that $\Lambda_{\eps_0}(q)<0$ for some $q\in C_0(\R)$ with $q\ge0$.
    Define
    \beaa
        \lambda
        :=\frac{\Lambda_{\eps_0}(\rho)}
        {-\Lambda_{\eps_0}(q)}>0,
        \qquad
        h:=-\rho-\lambda q.
    \eeaa
    Then $h<0$ and $\Lambda_{\eps_0}(h)=0$.
    Thus, $h\in\operatorname{Ran}L_{\eps_0}$, contradicting Lemma~\ref{lem:kernel-L}.
    Therefore, $\Lambda_{\eps_0}$ is a positive functional.

    Finally, Lemma~\ref{lem:Dalpha-phi-C0} gives $\dalpha\phi_{\eps_0}\in C_0(\R)$.
    Since $\phi_{\eps_0}'>0$, the definition of $\dalpha$ also gives $\dalpha\phi_{\eps_0}>0$.
    Applying Lemma~\ref{lem:kernel-L} to $-\dalpha\phi_{\eps_0}$, we obtain
    \beaa
        -\dalpha\phi_{\eps_0}
        \notin\operatorname{Ran}L_{\eps_0}
        =\ker\Lambda_{\eps_0}.
    \eeaa
    Hence $\Lambda_{\eps_0}(\dalpha\phi_{\eps_0})\ne0$.
    Together with the positivity of $\Lambda_{\eps_0}$, this yields
    $\Lambda_{\eps_0}(\dalpha\phi_{\eps_0})>0$.
\end{proof}

\medskip
We now prove Proposition~\ref{prop:extension}.
\begin{proof}[Proof of Proposition~\ref{prop:extension}]
    For $\sigma\in\R$ and $z=(\psi,\tau)\in C_0^2(\R)\times\R$, define
    \beaa
        H_{\eps_0}(\sigma,z)
        :=
        \left(
        (\phi_{\eps_0}+\psi)''
        -(\theta_{\eps_0}+\tau)
        \dalpha(\phi_{\eps_0}+\psi)
        +f_{\eps_0+\sigma}(\phi_{\eps_0}+\psi),
        \psi(0)
        \right).
    \eeaa
    Then $H_{\eps_0}:\R\times(C_0^2(\R)\times\R)\to C_0(\R)\times\R$ is of class $C^1$.
    Moreover, $H_{\eps_0}(0,0)=(0,0)$, and
    \beaa
        DH_{\eps_0}
        :=
        \frac{\partial H_{\eps_0}}{\partial z}(0,0)
        =
        \begin{bmatrix}
            L_{\eps_0} & -\dalpha\phi_{\eps_0}\\
            E_0 & 0
        \end{bmatrix}.
    \eeaa

    We show that $DH_{\eps_0}:C_0^2(\R)\times\R\to C_0(\R)\times\R$ is an isomorphism.
    For any $(A,b)\in C_0(\R)\times\R$, consider
    \beaa
        L_{\eps_0}\psi
        =A+\tau\dalpha\phi_{\eps_0},
        \qquad
        \psi(0)=b.
    \eeaa
    Since $\operatorname{Ran}L_{\eps_0}=\ker\Lambda_{\eps_0}$, the first equation has a solution $\psi\in C_0^2(\R)$ if and only if
    \beaa
        \Lambda_{\eps_0}(A)
        +\tau\Lambda_{\eps_0}
        \bigl(\dalpha\phi_{\eps_0}\bigr)=0.
    \eeaa
    By Lemma~\ref{lem:functional}, this condition uniquely determines $\tau$.
    By Lemma~\ref{lem:kernel-L}, any two solutions for this $\tau$ differ by a constant multiple of $\phi_{\eps_0}'$.
    Since $\phi_{\eps_0}'(0)>0$, the normalization $\psi(0)=b$ determines $\psi$ uniquely.
    Therefore, $DH_{\eps_0}$ is a bounded linear bijection.
    By the bounded inverse theorem, it is an isomorphism.

    The implicit function theorem therefore yields $\eta>0$ and a $C^1$ map
    \beaa
        (-\eta,\eta)\ni\sigma
        \longmapsto
        (\psi_\sigma,\tau_\sigma)\in C_0^2(\R)\times\R
    \eeaa
    satisfying $(\psi_0,\tau_0)=(0,0)$ and
    $H_{\eps_0}(\sigma,(\psi_\sigma,\tau_\sigma))=(0,0)$.
    Reducing $\eta$ if necessary, we may assume that $\eta<1-\eps_0$.
    Then $f_\eps$ satisfies \eqref{ass:nonlinearity} for all $\eps\in[\eps_0,\eps_0+\eta)$.
    For these $\eps$, define
    \beaa
        \phi_\eps
        :=\phi_{\eps_0}+\psi_{\eps-\eps_0},
        \qquad
        \theta_\eps
        :=\theta_{\eps_0}+\tau_{\eps-\eps_0}.
    \eeaa
    The profile equation and normalization in Problem~\eqref{eq:problem-atw} follow from $H_{\eps_0}(\eps, (\psi_{\eps-\eps_0}, \tau_{\eps-\eps_0}))=(0,0)$.
    The conditions \eqref{ass:boundary_condition} follow from $\psi_{\eps-\eps_0}\in C_0^2(\R)$.
    Proposition~\ref{prop:posi-speed} gives $\theta_\eps>0$.
    Thus, $(\theta_\eps,\phi_\eps)$ is a $C^1$ family of solutions to Problem~\eqref{eq:problem-atw}.
\end{proof}

\medskip
\subsection{Proof of the main theorem}\label{subsec:proof}
We now prove Theorem~\ref{thm:main}.

Proposition~\ref{prop:local-zero} gives a local family of solutions starting from $(0,\phi_0)$ at $\eps=0$.
Define
\beaa
    \eps_*
    :=
    \sup
    \left\{
    \eta\in(0,1] \;\middle|\;
    \begin{array}{l}
    \text{for every $\eps\in[0,\eta)$, Problem \eqref{eq:problem-atw} has a solution}\\
    (\theta_\eps,\phi_\eps)\in\R\times C_b^2(\R).
    \end{array}
    \right\}.
\eeaa
By Proposition~\ref{prop:local-zero}, we have $0<\eps_*\le1$.
Moreover, Propositions~\ref{prop:posi-speed} and~\ref{prop:mono} give $\phi_\eps'>0$ for $\eps>0$, while $\phi_0'>0$.
By the definition of $\eps_*$, solutions exist for every $\eps\in[0,\eps_*)$.
Proposition~\ref{prop:limit-solution} also gives a solution at $\eps=\eps_*$.

Suppose that $\eps_*<1$.
Applying Proposition~\ref{prop:extension} with $\eps_0=\eps_*$, we can extend the interval of existence.
More precisely, there exists $\eta\in(0,1-\eps_*)$ such that, for every $\eps\in(0,\eps_*+\eta)$, Problem~\eqref{eq:problem-atw} has a solution $(\theta_\eps,\phi_\eps)\in(0,+\infty)\times C_b^2(\R)$.
This contradicts the definition of $\eps_*$.
Therefore, $\eps_*=1$.

Since $f_1=f$, this yields the required asymptotic traveling wave with $c=\theta_1^{1/\alpha}>0$.
Strict monotonicity follows from Proposition~\ref{prop:mono}.
Uniqueness up to translation follows from Proposition~\ref{prop:uni} and the relation $c=\theta^{1/\alpha}$.
This completes the proof of Theorem~\ref{thm:main}.

\bigskip
\section{Speed asymptotics and numerical simulations}\label{sec:num}

In this section, we consider Equation~\eqref{eq:problem} with a cubic nonlinearity satisfying \eqref{ass:nonlinearity}.
We examine the relationship between the asymptotic wave speed and the numerically computed front speed.
We first obtain an explicit asymptotic speed formula using Remark~\ref{rem:speed-zero}.
We then describe the numerical setting and compare the computed front speed with the leading term of this formula.

\subsection{Asymptotic speed near the balanced case}\label{subsec:sp-atw}

In this subsection, we consider Equation~\eqref{eq:problem} with the nonlinearity
\beaa
    f_\delta(u):=u(1-u)\left(u-\frac{1}{2}+\delta\right),
    \qquad 0\le\delta<\frac{1}{2}.
\eeaa
We study the asymptotic wave speed near the balanced case.
This family differs from $f_\eps$ constructed in the preceding sections.
For $0<\delta<1/2$, $f_\delta$ satisfies \eqref{ass:nonlinearity} with $a=1/2-\delta$.
The case $\delta=0$ is balanced.

In the classical case $\alpha=1$, the equation $u_t=u_{xx}+f_\delta(u)$ admits the explicit Huxley traveling front
\beaa
    u(t,x)=\Phi_0(x+c_\delta^{\mathrm{cl}}t),
    \qquad \Phi_0(\xi)=\frac{1}{1+\exp(-\xi/\sqrt2)},
    \qquad
    c_\delta^{\mathrm{cl}}=\sqrt2\,\delta.
\eeaa
At $\delta=0$, $\Phi_0$ is also a stationary front of the time-fractional equation and satisfies $\Phi_0(0)=1/2$.
We now fix $\alpha\in(0,1)$ and use $\Phi_0$ as the reference profile for the following asymptotic analysis.

Let $\Phi_\delta$ and $c_\delta$ denote the profile and speed of the asymptotic traveling wave corresponding to $f_\delta$.
We impose the normalization $\Phi_\delta(0)=1/2$ and set $\theta_\delta:=c_\delta^\alpha$.
As in Proposition~\ref{prop:local-zero}, we apply the implicit function theorem at $(\delta,\Phi,\theta)=(0,\Phi_0,0)$.
This gives a family of solutions for sufficiently small $\delta\ge0$.
Moreover, the map
$\delta\mapsto(\Phi_\delta-\Phi_0,\theta_\delta)\in C_0^2(\R)\times\R$
is of class $C^1$.

As in Remark~\ref{rem:speed-zero}, we differentiate the profile equation at $\delta=0$ and apply the solvability condition.
Since $\partial_\delta f_\delta(u)=u(1-u)$, we obtain
\beaa
    \left.\frac{d\theta_\delta}{d\delta}\right|_{\delta=0}
    =
    \frac{\dint_{\R}\Phi_0(1-\Phi_0)\Phi_0'\,d\xi}
    {\dint_{\R}(\dalpha\Phi_0)\Phi_0'\,d\xi}
    =
    \frac{1}{6I_\alpha}
\eeaa
where $I_\alpha:=\int_{\R}(\dalpha\Phi_0)\Phi_0'\,d\xi$.

We now compute $I_\alpha$ explicitly.
We use the Fourier transform convention
\beaa
    \widehat{f}(k)
    :=
    \mathcal{F}[f](k)
    =
    \frac{1}{\sqrt{2\pi}}
    \int_{\R}f(x)e^{-ikx}\,dx.
\eeaa
Note that $\Phi_0'$ belongs to the Schwartz space $\cS(\R)$.
By the definition of $\dalpha$,
\beaa
    \dalpha\Phi_0=K_\alpha*\Phi_0',
    \qquad
    K_\alpha(s):=
    \frac{s^{-\alpha}}{\Gamma(1-\alpha)}
    \mathbf{1}_{(0,\infty)}(s).
\eeaa
It holds that $\widehat{K_\alpha}(k)=(2\pi)^{-1/2}(ik)^{\alpha-1}$ in the sense of tempered distributions, where complex powers are defined using the principal branch of the logarithm.
Hence
\be\label{eq:fourier}
    \mathcal F[\dalpha\Phi_0](k)
    =(ik)^{\alpha-1}\widehat{\Phi_0'}(k).
\ee
Since $\Phi_0'\in\mathcal S(\R)$, Parseval's identity in the sense of tempered distributions gives
\be\label{eq:parseval}
    I_\alpha
    =
    \int_{\R}
    \mathcal F[\dalpha\Phi_0](k)
    \overline{\widehat{\Phi_0'}(k)}\,dk
    =
    \int_{\R}(ik)^{\alpha-1}|\widehat{\Phi_0'}(k)|^2\,dk.
\ee
A direct calculation gives
\beaa
    \widehat{\Phi_0'}(k)
    =
    \frac{\sqrt{\pi}\,k}{\sinh(\sqrt{2}\pi k)}.
\eeaa
Using these formulas, we obtain the following result.
\begin{proposition}\label{prop:theta_alpha}
    The asymptotic wave speed $c_\delta$ satisfies
    \be\label{eq:num-speed-asy}
        c_\delta
        =C_\alpha\delta^{1/\alpha}+o(\delta^{1/\alpha})
        \qquad(\delta \downarrow 0),
    \ee
    where
    \be\label{eq:theta}
        C_\alpha
        :=\left(
            \frac{2^{3\alpha/2}\pi^{\alpha+1}}
            {6\sin(\pi\alpha/2)\Gamma(\alpha+2)\zeta(\alpha+1)}
        \right)^{1/\alpha}
    \ee
    and $\zeta$ is the Riemann zeta function.
\end{proposition}
\begin{proof}
    Since $\widehat{\Phi_0'}$ is real-valued and even,
    \eqref{eq:fourier} and \eqref{eq:parseval} give
    \beaa
        I_\alpha
        =
        \int_{\R}(ik)^{\alpha-1}
        |\widehat{\Phi_0'}(k)|^2\,dk
        &=&
        2\pi\sin\left(\frac{\pi\alpha}{2}\right)
        \int_0^\infty
        \frac{k^{\alpha+1}}{\sinh^2(\sqrt{2}\pi k)}\,dk\\
        &=&
        \frac{2\pi\sin(\pi\alpha/2)}
             {(\sqrt{2}\pi)^{\alpha+2}}
        \int_0^\infty\frac{x^{\alpha+1}}{\sinh^2x}\,dx.
    \eeaa
    For $x>0$, we have
    \beaa
        \frac{1}{\sinh^2x}
        =
        \frac{4e^{-2x}}{(1-e^{-2x})^2}
        =
        4\sum_{n=1}^\infty n e^{-2nx}.
    \eeaa
    Since all terms are nonnegative, the monotone convergence theorem and the integral representation of the gamma function yield
    \beaa
        \int_0^\infty\frac{x^{\alpha+1}}{\sinh^2x}\,dx
        &=&
        4\sum_{n=1}^\infty n
        \int_0^\infty x^{\alpha+1}e^{-2nx}\,dx\\    
        &=&
        4\sum_{n=1}^\infty
        \frac{n\Gamma(\alpha+2)}{(2n)^{\alpha+2}}
        =
        2^{-\alpha}\Gamma(\alpha+2)\zeta(\alpha+1).
    \eeaa
    Simplifying the above expressions, 
    we find that $(6I_\alpha)^{-1/\alpha}$ equals $C_\alpha$ defined in \eqref{eq:theta}.
\end{proof}
\begin{remark}
    The right-hand side of \eqref{eq:theta} is defined for $\alpha\in(0,1)$.
    Letting $\alpha\uparrow 1$, we obtain
    \beaa
        \lim_{\alpha\uparrow1 }C_\alpha
        =
        \frac{2^{\frac{3}{2}}\pi^2}
        {6\cdot1\cdot2\cdot\pi^2/6}
        =
        \sqrt{2}.
    \eeaa
    This agrees with the coefficient of $\delta$ in the exact speed formula $c_\delta^{\mathrm{cl}}=\sqrt2\,\delta$.
\end{remark}

\begin{figure}[tb]
  \centering
  \includegraphics[width=0.65\textwidth]{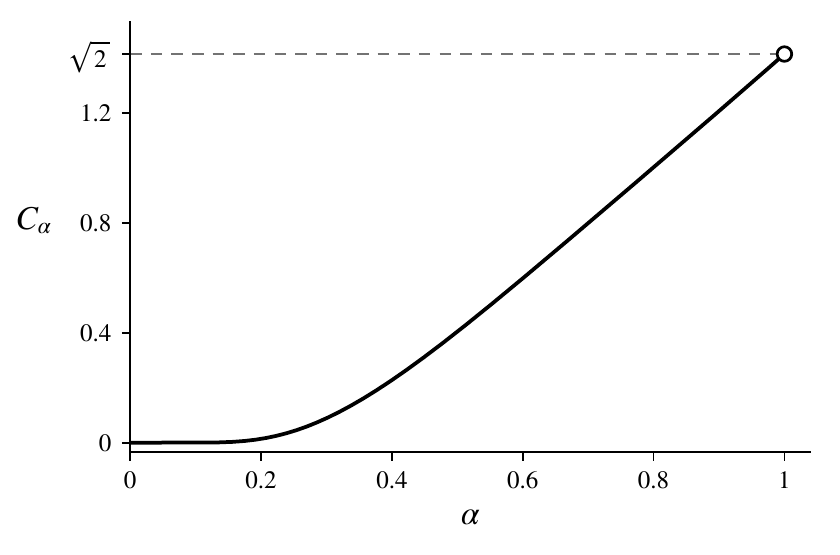}
  \caption{\small{The coefficient $C_{\alpha}$ in
  \eqref{eq:num-speed-asy} as a function of $\alpha$.
  The dashed line indicates the limiting value $\sqrt{2}$
  as $\alpha\uparrow1$, which is marked by the open circle.}}
  \label{fig:coefficient}
\end{figure}

\medskip
For the cubic nonlinearity, we obtain the asymptotic speed formula \eqref{eq:num-speed-asy}.
Figure~\ref{fig:coefficient} shows $C_{\alpha}$ as a function of $\alpha$.
The plot suggests that $C_{\alpha}$ is strictly increasing in $\alpha$.

\medskip
\subsection{Numerical scheme and results}\label{subsec:res}

We describe the numerical setting and results.
The asymptotic formula \eqref{eq:num-speed-asy} gives $c_\delta=O(\delta^{1/\alpha})$ for each fixed $\alpha\in(0,1)$.
To observe front propagation near the balanced case, we rescale the equation as $t\mapsto \delta^{-(1+\alpha)/\alpha}t$ and $x\mapsto \delta^{-1} x$ so that the asymptotic wave speed is $O(1)$.
For the simulation, we solve the problem
\beaa
    \begin{cases}
        \delta^{1+\alpha}\partial_t^\alpha u
        =\delta^2u_{xx}+u(1-u)\left(u-\frac12+\delta\right)
        \quad(t>0,\ 0<x<l),\\
        u_x(t,0)=u_x(t,l)=0
        \quad(t>0),\\
        u(0,x)=u_0(x)
        \quad(0<x<l),
    \end{cases}
\eeaa
where $l>0$ and 
\beaa
    u_0(x)=\frac12\left(1+\tanh\left(\frac{x-l_0}{\delta}\right)\right)
\eeaa
for some $l_0 \in (0,l)$.
We use the numerical scheme described in \cite[Section~4.1]{Ishii2025}.
This scheme combines the L1 approximation of the Caputo derivative with a central difference approximation of the second spatial derivative.
The diffusion term is treated implicitly and the reaction term explicitly.
We use a spatial mesh size $\Delta x=6.0 \times 10^{-3}$ and a time step $\Delta t=4.0\times10^{-4}$.

In the rescaled coordinates, the asymptotic wave speed is $\delta^{-1/\alpha}c_\delta$.
For fixed $\alpha$, the asymptotic formula \eqref{eq:num-speed-asy} shows that its leading term as $\delta\downarrow0$ is $C_\alpha$.
Therefore, we compare the propagation speed of a numerical solution with $C_\alpha$.

We write $c_{\alpha,\delta}^{\mathrm{num}}$ for the numerical speed.
Let $T>0$ be the final simulation time.
We compute $c_{\alpha,\delta}^{\mathrm{num}}$ as follows:
\begin{enumerate}
    \item
    Compute the interface position $x^*(t)$ satisfying
    \beaa
        u(t,x^*(t))=1/2, \qquad (t>0)
    \eeaa
    by linear interpolation between adjacent grid points whose solution values bracket $1/2$.
    Record the position at time intervals of $0.2$.

    \item
    Fit the line $x^*(t)\simeq b+st$ to the data for $T/2\le t\le T$ using least squares.
    Set $c_{\alpha,\delta}^{\mathrm{num}} := |s|$.
\end{enumerate}

\begin{figure}[!tb]
  \centering
    \includegraphics[width=1.0\linewidth]{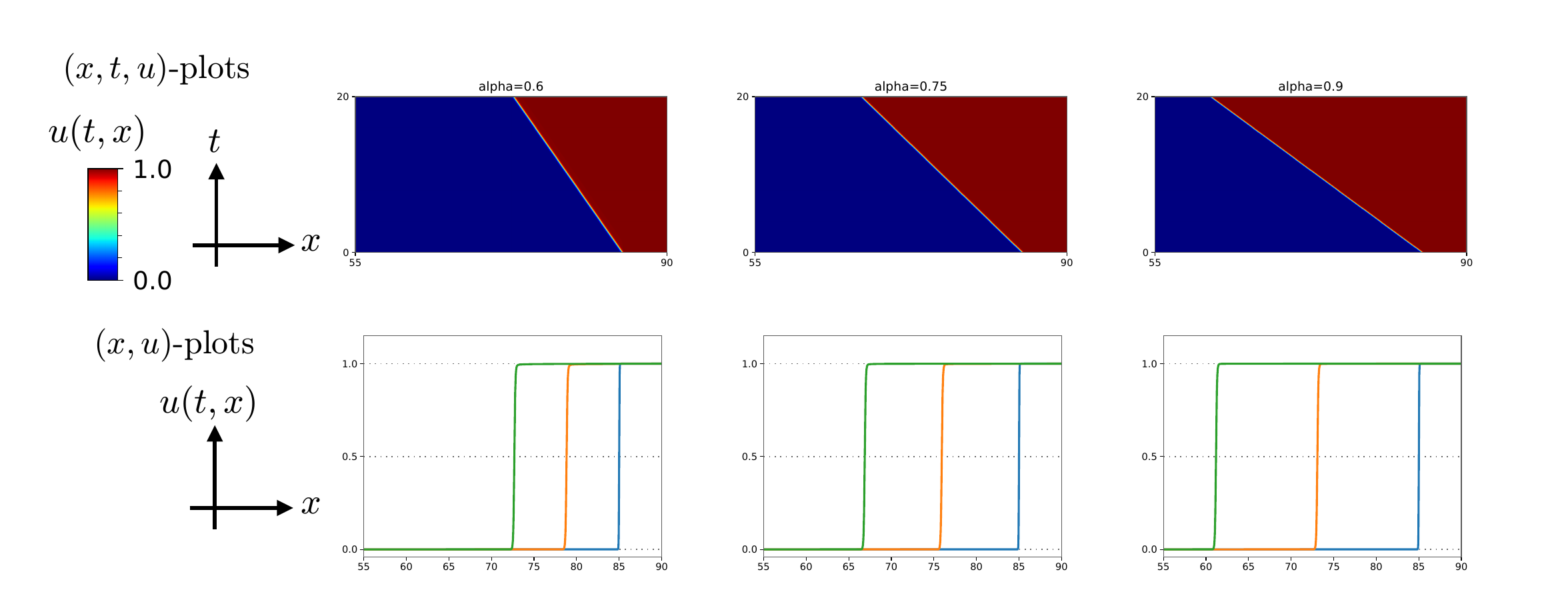}
  \caption{\small{
    Spatio-temporal plots (top panels) and profiles of the numerical solutions (bottom panels).
    The parameters are $l=90$, $l_0=85$, $\delta=0.03$, and $T=20$.
    The numerical simulations correspond to $\alpha=0.6$ (left), $\alpha=0.75$ (middle), and $\alpha=0.9$ (right).
    The profiles are plotted in blue ($t=0$), orange ($t=10$), and green ($t=20$).
  }
  }
  \label{fig:phase}
\end{figure}

The profile of the numerical solution is shown in Fig.~\ref{fig:phase}.
The spatio-temporal plots show a front connecting $u=0$ and $u=1$ moving from right to left.
The profiles at different times also show the leftward motion of the transition layer.
Over the time interval used for speed estimation, the front position varies approximately linearly with time.
This indicates propagation at an approximately constant speed over this interval.
Moreover, the front travels farther over the same time interval for larger $\alpha$.

\begin{figure}[!tb]
  \centering
    \includegraphics[width=1.0\linewidth]{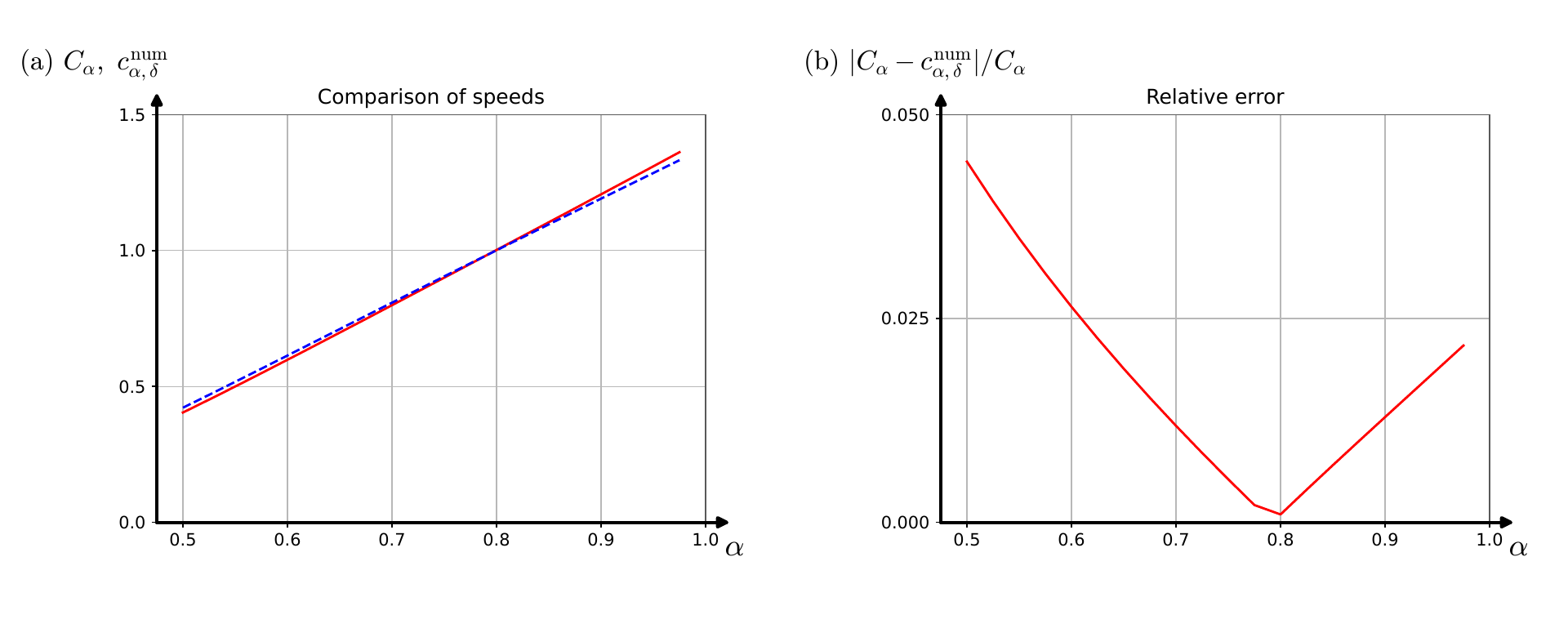}
  \caption{\small{
    (a) Comparison of $C_\alpha$ (red line) and the propagation speed of the numerical solution (dashed line).
    The parameters are $l=90$, $l_0=85$, $\delta=0.04$, and $T=20$.
    The  blue dashed line connects the numerical speeds computed  for 
    $\alpha = 0.5 + 0.025 \times k\, (0 \le k \le 19)$ by straight line segments.
    (b) Relative errors.
  }
  }
  \label{fig:speed-fixed-delta}
\end{figure}

We next compare the numerical speed $c_{\alpha,\delta}^{\mathrm{num}}$ with the leading term $C_\alpha$ of the rescaled speed.
Figure~\ref{fig:speed-fixed-delta} shows the results for varying $\alpha$ at fixed $\delta$.
Both $c_{\alpha,\delta}^{\mathrm{num}}$ and $C_\alpha$ increase with $\alpha$.
The relative error is less than $0.045$ and smallest near $\alpha=0.8$.

\begin{figure}[!tb]
  \centering
    \includegraphics[width=1.0\linewidth]{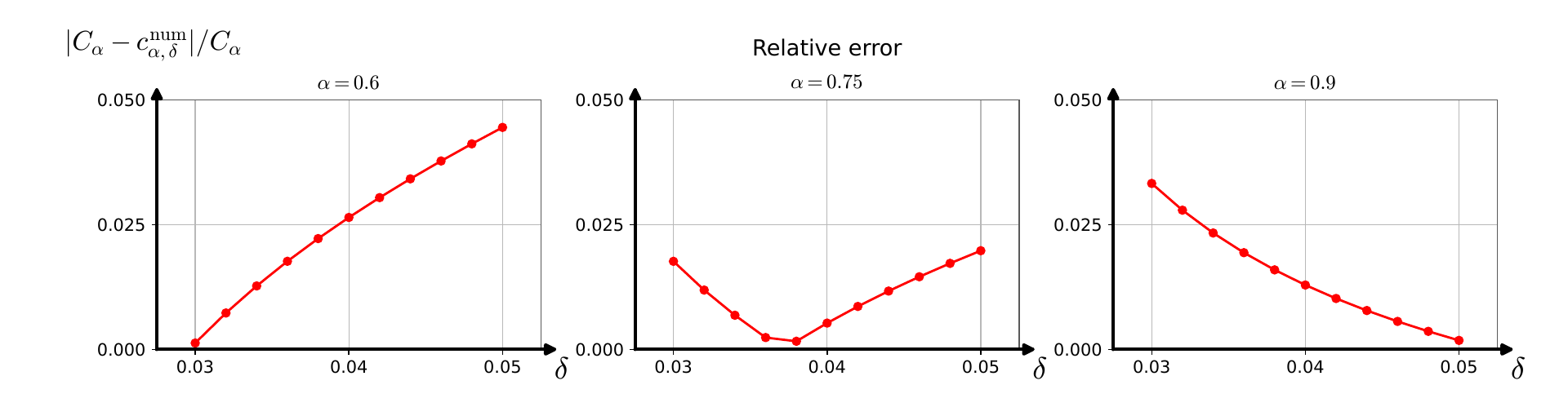}
    \caption{\small{
    Relative errors for varying $\delta$.
    The parameters are $l=90$, $l_0=85$, and $T=20$.
    The panels correspond to $\alpha=0.6$ (left),
    $\alpha=0.75$ (middle), and $\alpha=0.9$ (right).
    }}
  \label{fig:speed-fixed-alpha}
\end{figure}

The asymptotic formula \eqref{eq:num-speed-asy} holds as $\delta\downarrow0$ for each fixed $\alpha$.
Figure~\ref{fig:speed-fixed-alpha} shows the relative errors for varying $\delta$ at fixed $\alpha$.
All relative errors are below $0.05$, but they do not consistently decrease as $\delta$ decreases.

These results indicate that $c^{\mathrm{num}}_{\alpha,\delta}$ is close to $C_\alpha$ over the parameter ranges and time intervals considered.
Since $C_\alpha$ is the leading-order coefficient in \eqref{eq:num-speed-asy}, the observed differences may reflect finite-$\delta$ effects, finite-time effects, and numerical discretization.
Some computations at smaller $\delta$ exhibited instability with the discretization used here; we thus restrict the comparison to the ranges shown in the figures.
Further investigation is needed to assess this correspondence as $\delta\downarrow0$.

\bigskip
\section{Concluding remarks}\label{sec:dis}
In this paper, we analyzed the asymptotic profile equation derived from the Caputo time-fractional reaction-diffusion equation.
For general $C^1$ bistable nonlinearities satisfying \eqref{ass:nonlinearity}, we proved the existence and uniqueness up to translation of asymptotic traveling waves.
We also established strict monotonicity of the profiles and positivity of the wave speed under the imbalance condition.
The existence proof used a continuation method starting from a stationary front in the balanced case.

For a specific cubic nonlinearity, we also derived the asymptotic speed formula \eqref{eq:num-speed-asy} near the balanced case.
We obtained an explicit expression for its coefficient $C_\alpha$.
For each fixed $\alpha\in(0,1)$, the speed satisfies
$c_\delta\sim C_\alpha\delta^{1/\alpha}\ (\delta\downarrow0)$.
Thus, the order $\alpha$ of the Caputo derivative determines
the scaling exponent $1/\alpha$.
Numerical simulations also showed front propagation at an approximately constant speed over the observed time interval.
The rescaled numerical speed $c_{\alpha,\delta}^{\mathrm{num}}$ was close to $C_\alpha$ over the range considered.
In comparisons at fixed $\delta$, both $c_{\alpha,\delta}^{\mathrm{num}}$ and $C_\alpha$ increased with $\alpha$.

Further work is needed to analyze the long-time behavior of solutions to the Cauchy problem for Equation~\eqref{eq:problem}.
For the time-fractional Amari model, solutions of the Cauchy problem with suitable initial data converge to an asymptotic traveling wave profile.
This convergence has been proved in a moving coordinate centered at the front position~\cite{Ishii-preprint}.
Whether analogous convergence holds for Equation~\eqref{eq:problem} remains open.
However, as noted in Section~\ref{sec:intro}, asymptotic traveling waves are not solutions of Equation~\eqref{eq:problem}.
This makes it difficult to apply stability methods for traveling waves in classical reaction-diffusion equations directly.
In the Amari model, a Heaviside-type nonlinearity allows a reduction to a Volterra equation for the front position.
The convergence proof relies on this reduction, which is difficult to carry out for the bistable nonlinearities considered here.
Establishing convergence to asymptotic traveling waves for suitable initial data and determining the convergence rate remain important open problems.

\bigskip
\section*{Acknowledgments}
HI is supported by JSPS KAKENHI Grant Numbers JP23K13013, JP24H00188 and JP25K07142.

\bigskip
\section*{Declaration of generative AI and AI-assisted technologies in the manuscript preparation process}
During the preparation of this work, the authors used OpenAI's ChatGPT and Codex, as well as Google's Gemini to assist with literature searches, exploration of possible proof strategies, examination of mathematical arguments, English translation, language editing, and the preparation of scripts for managing numerical simulations and generating figures. After using these tools, the authors reviewed and edited the generated content as needed and take full responsibility for the content of the publication.

\bigskip
\section*{Declaration of competing interests}
The authors declare no competing interests.

\bibliographystyle{abbrv}
\bibliography{export}

@article{Nec2010,
  author        = {Nec, Yana and Volpert, Vladimir A. and Nepomnyashchy, Alexander A.},
  title         = {Front propagation problems with sub-diffusion},
  journal       = {Discrete Contin. Dyn. Syst.},
  volume        = {27},
  number        = {2},
  year          = {2010},
  pages         = {827--846},
  doi           = {10.3934/dcds.2010.27.827},
}

@article{Bates1997,
  author        = {Bates, Peter W. and Fife, Paul C. and Ren, Xiaofeng and Wang, Xuefeng},
  title         = {Traveling waves in a convolution model for phase transitions},
  journal       = {Arch. Ration. Mech. Anal.},
  volume        = {138},
  number        = {2},
  year          = {1997},
  pages         = {105--136},
  doi           = {10.1007/s002050050037},
}

@article{Ishii2025,
  author        = {Ishii, Hiroshi},
  title         = {Propagating front solutions in a time-fractional {Fisher--KPP} equation},
  journal       = {Discrete Contin. Dyn. Syst. Ser. B},
  volume        = {30},
  number        = {7},
  year          = {2025},
  pages         = {2460--2482},
  doi           = {10.3934/dcdsb.2024173},
}

@article{Alfaro2008,
  author        = {Alfaro, Matthieu and Hilhorst, Danielle and Matano, Hiroshi},
  title         = {The singular limit of the {Allen--Cahn} equation and the {FitzHugh--Nagumo} system},
  journal       = {J. Differential Equations},
  volume        = {245},
  number        = {2},
  year          = {2008},
  pages         = {505--565},
  doi           = {10.1016/j.jde.2008.01.014},
}

@article{Bates2006,
  author        = {Bates, Peter W. and Chen, Fengxin},
  title         = {Spectral analysis of traveling waves for nonlocal evolution equations},
  journal       = {SIAM J. Math. Anal.},
  volume        = {38},
  number        = {1},
  year          = {2006},
  pages         = {116--126},
  doi           = {10.1137/S0036141004443968},
}

@incollection{Hille1972,
  author        = {Hille, Einar},
  title         = {Generalizations of {Landau's} inequality to linear operators},
  editor        = {Butzer, P. L. and Kahane, J.-P. and Sz{\"o}kefalvi-Nagy, B.},
  booktitle     = {Linear Operators and Approximation / Lineare Operatoren und Approximation},
  series        = {International Series of Numerical Mathematics},
  volume        = {20},
  publisher     = {Birkh{\"a}user},
  address       = {Basel},
  year          = {1972},
  pages         = {20--32},
  doi           = {10.1007/978-3-0348-7283-6_1},
}

@book{GT,
  author        = {Gilbarg, David and Trudinger, Neil S.},
  title         = {Elliptic Partial Differential Equations of Second Order},
  series        = {Classics in Mathematics},
  edition       = {Second},
  publisher     = {Springer},
  address       = {Berlin, Heidelberg},
  year          = {2001},
  doi           = {10.1007/978-3-642-61798-0},
}

@book{Volpert1994,
  author        = {Volpert, Aizik I. and Volpert, Vitaly A. and Volpert, Vladimir A.},
  title         = {Traveling Wave Solutions of Parabolic Systems},
  series        = {Translations of Mathematical Monographs},
  volume        = {140},
  publisher     = {American Mathematical Society},
  address       = {Providence, RI},
  year          = {1994},
  doi           = {10.1090/mmono/140},
}

@article{QTWY,
  author        = {Quan, Chaoyu and Tang, Tao and Wang, Boyi and Yang, Jiang},
  title         = {A decreasing upper bound of the energy for time-fractional phase-field equations},
  journal       = {Commun. Comput. Phys.},
  volume        = {33},
  number        = {4},
  year          = {2023},
  pages         = {962--991},
  doi           = {10.4208/cicp.OA-2022-0148},
}

@article{AC97,
  author        = {Allen, Samuel M. and Cahn, John W.},
  title         = {A microscopic theory for antiphase boundary motion and its application to antiphase domain coarsening},
  journal       = {Acta Metall.},
  volume        = {27},
  number        = {6},
  year          = {1979},
  pages         = {1085--1095},
  doi           = {10.1016/0001-6160(79)90196-2},
}

@article{NYA65,
  author        = {Nagumo, J. and Yoshizawa, S. and Arimoto, S.},
  title         = {Bistable transmission lines},
  journal       = {IEEE Trans. Circuit Theory},
  volume        = {12},
  number        = {3},
  year          = {1965},
  pages         = {400--412},
  doi           = {10.1109/TCT.1965.1082476},
}

@article{CH58,
  author        = {Cahn, John W. and Hilliard, John E.},
  title         = {Free energy of a nonuniform system. {I}. {Interfacial} free energy},
  journal       = {J. Chem. Phys.},
  volume        = {28},
  number        = {2},
  year          = {1958},
  pages         = {258--267},
  doi           = {10.1063/1.1744102},
}

@article{FM77,
  author        = {Fife, Paul C. and McLeod, John Bryce},
  title         = {The approach of solutions of nonlinear diffusion equations to travelling front solutions},
  journal       = {Arch. Ration. Mech. Anal.},
  volume        = {65},
  number        = {4},
  year          = {1977},
  pages         = {335--361},
  doi           = {10.1007/BF00250432},
  url           = {https://api.semanticscholar.org/CorpusID:10307450},
}

@article{FMN04,
  author        = {Fukao, Yukitaka and Morita, Yoshihisa and Ninomiya, Hirokazu},
  title         = {Some entire solutions of the {Allen--Cahn} equation},
  journal       = {Taiwanese J. Math.},
  volume        = {8},
  number        = {1},
  year          = {2004},
  pages         = {15--32},
  doi           = {10.11650/twjm/1500558454},
}

@article{JM05,
  author        = {Guo, Jong-Shenq and Morita, Yoshihisa},
  title         = {Entire solutions of reaction-diffusion equations and an application to discrete diffusive equations},
  journal       = {Discrete Contin. Dyn. Syst.},
  volume        = {12},
  number        = {2},
  year          = {2005},
  pages         = {193--212},
  doi           = {10.3934/dcds.2005.12.193},
}

@article{CGN06,
  author        = {Chen, Xinfu and Guo, Jong-Shenq and Ninomiya, Hirokazu},
  title         = {Entire solutions of reaction---diffusion equations with balanced bistable nonlinearities},
  journal       = {Proc. Roy. Soc. Edinburgh Sect. A},
  volume        = {136},
  number        = {6},
  year          = {2006},
  pages         = {1207--1237},
  doi           = {10.1017/S0308210500004959},
}

@article{MN06,
  author        = {Morita, Yoshihisa and Ninomiya, Hirokazu},
  title         = {Entire solutions with merging fronts to reaction--diffusion equations},
  journal       = {J. Dynam. Differential Equations},
  volume        = {18},
  number        = {4},
  year          = {2006},
  pages         = {841--861},
  doi           = {10.1007/s10884-006-9046-x},
}

@article{HMR05,
  author        = {Hamel, Fran{\c{c}}ois and Monneau, R{\'e}gis and Roquejoffre, Jean-Michel},
  title         = {Existence and qualitative properties of multidimensional conical bistable fronts},
  journal       = {Discrete Contin. Dyn. Syst.},
  volume        = {13},
  number        = {4},
  year          = {2005},
  pages         = {1069--1096},
  doi           = {10.3934/dcds.2005.13.1069},
}

@article{HMR06,
  author        = {Hamel, Fran{\c{c}}ois and Monneau, R{\'e}gis and Roquejoffre, Jean-Michel},
  title         = {Asymptotic properties and classification of bistable fronts with {Lipschitz} level sets},
  journal       = {Discrete Contin. Dyn. Syst.},
  volume        = {14},
  number        = {1},
  year          = {2006},
  pages         = {75--92},
  doi           = {10.3934/dcds.2006.14.75},
}

@article{NT05,
  author        = {Ninomiya, Hirokazu and Taniguchi, Masaharu},
  title         = {Existence and global stability of traveling curved fronts in the {Allen--Cahn} equations},
  journal       = {J. Differential Equations},
  volume        = {213},
  number        = {1},
  year          = {2005},
  pages         = {204--233},
  doi           = {10.1016/j.jde.2004.06.011},
}

@article{NT06,
  author        = {Ninomiya, Hirokazu and Taniguchi, Masaharu},
  title         = {Global stability of traveling curved fronts in the {Allen-Cahn} equations},
  journal       = {Discrete Contin. Dyn. Syst.},
  volume        = {15},
  number        = {3},
  year          = {2006},
  pages         = {819--832},
  doi           = {10.3934/dcds.2006.15.819},
}

@article{BSS93,
  author        = {Barles, G. and Soner, H. M. and Souganidis, P. E.},
  title         = {Front propagation and phase field theory},
  journal       = {SIAM J. Control Optim.},
  volume        = {31},
  number        = {2},
  year          = {1993},
  pages         = {439--469},
  doi           = {10.1137/0331021},
}

@article{BBS92,
  author        = {Barles, G. and Bronsard, L. and Souganidis, P. E.},
  title         = {Front propagation for reaction-diffusion equations of bistable type},
  journal       = {Ann. Inst. H. Poincar{\'e} Anal. Non Lin{\'e}aire},
  volume        = {9},
  number        = {5},
  year          = {1992},
  pages         = {479--496},
  doi           = {10.1016/S0294-1449(16)30228-1},
}

@article{MK00,
  author        = {Metzler, Ralf and Klafter, Joseph},
  title         = {The random walk's guide to anomalous diffusion: a fractional dynamics approach},
  journal       = {Phys. Rep.},
  volume        = {339},
  number        = {1},
  year          = {2000},
  pages         = {1--77},
  doi           = {10.1016/S0370-1573(00)00070-3},
}

@article{HWL06,
  author        = {Henry, B. I. and Langlands, T. A. M. and Wearne, S. L.},
  title         = {Anomalous diffusion with linear reaction dynamics: from continuous time random walks to fractional reaction-diffusion equations},
  journal       = {Phys. Rev. E (3)},
  volume        = {74},
  number        = {3},
  year          = {2006},
  note          = {Article 031116},
  doi           = {10.1103/PhysRevE.74.031116},
}

@article{FK26,
  author        = {Fritz, Marvin and Kavallaris, Nikos I.},
  title         = {A time-fractional {Fisher--KPP} equation for tumor growth: analysis and numerical simulation},
  journal       = {Commun. Nonlinear Sci. Numer. Simul.},
  volume        = {159},
  year          = {2026},
  note          = {Article 109911},
  doi           = {10.1016/j.cnsns.2026.109911},
}

@article{AH20,
  author        = {Angstmann, Christopher N. and Henry, Bruce I.},
  title         = {Time fractional {Fisher--KPP} and {Fitzhugh--Nagumo} equations},
  journal       = {Entropy},
  volume        = {22},
  number        = {9},
  year          = {2020},
  note          = {Article 1035},
  doi           = {10.3390/e22091035},
  url           = {https://www.mdpi.com/1099-4300/22/9/1035},
}

@article{VNN13,
  author        = {Volpert, V. A. and Nec, Y. and Nepomnyashchy, Alexander A.},
  title         = {Fronts in anomalous diffusion--reaction systems},
  journal       = {Philos. Trans. Roy. Soc. A},
  volume        = {371},
  number        = {1982},
  year          = {2013},
  note          = {Article 20120179},
  doi           = {10.1098/rsta.2012.0179},
}

@misc{Ishii-preprint,
  author        = {Ishii, Hiroshi},
  title         = {Global front solutions converging to an asymptotic traveling wave in a time-fractional {Amari} model},
  year          = {2026},
  howpublished  = {arXiv:2609.14884 [math.AP]},
  doi           = {10.48550/arXiv.2609.14884},
  eprint        = {2609.14884},
  archiveprefix = {arXiv},
  primaryclass  = {math.AP},
  url           = {https://arxiv.org/abs/2609.14884},
}

@article{Gonzalez,
  author        = {Gonz{\'a}lez-Ram{\'i}rez, L. R.},
  title         = {On the existence of traveling fronts in the fractional-order {Amari} neural field model},
  journal       = {Commun. Nonlinear Sci. Numer. Simul.},
  volume        = {116},
  year          = {2023},
  note          = {Article 106790},
  doi           = {10.1016/j.cnsns.2022.106790},
}

@article{LCWZ18,
  author        = {Liu, Huan and Cheng, Aijie and Wang, Hong and Zhao, Jia},
  title         = {Time-fractional {Allen--Cahn} and {Cahn--Hilliard} phase-field models and their numerical investigation},
  journal       = {Comput. Math. Appl.},
  volume        = {76},
  number        = {8},
  year          = {2018},
  pages         = {1876--1892},
  doi           = {10.1016/j.camwa.2018.07.036},
}

@article{DYZ20,
  author        = {Du, Qiang and Yang, Jiang and Zhou, Zhi},
  title         = {Time-fractional {Allen--Cahn} equations: analysis and numerical methods},
  journal       = {J. Sci. Comput.},
  volume        = {85},
  number        = {2},
  year          = {2020},
  note          = {Article 42},
  doi           = {10.1007/s10915-020-01351-5},
}

@article{TYZ19,
  author        = {Tang, Tao and Yu, Haijun and Zhou, Tao},
  title         = {On energy dissipation theory and numerical stability for time-fractional phase-field equations},
  journal       = {SIAM J. Sci. Comput.},
  volume        = {41},
  number        = {6},
  year          = {2019},
  pages         = {A3757--A3778},
  doi           = {10.1137/18M1203560},
}

@article{EM93,
  author        = {Ermentrout, G. Bard and McLeod, J. Bryce},
  title         = {Existence and uniqueness of travelling waves for a neural network},
  journal       = {Proc. Roy. Soc. Edinburgh Sect. A},
  volume        = {123},
  number        = {3},
  year          = {1993},
  pages         = {461--478},
  doi           = {10.1017/S030821050002583X},
}

@article{GUI15,
  author        = {Gui, Changfeng and Zhao, Mingfeng},
  title         = {Traveling wave solutions of {Allen--Cahn} equation with a fractional {Laplacian}},
  journal       = {Ann. Inst. H. Poincar{\'e} Anal. Non Lin{\'e}aire},
  volume        = {32},
  number        = {4},
  year          = {2015},
  pages         = {785--812},
  doi           = {10.1016/j.anihpc.2014.03.005},
}

@article{Chen97,
  author        = {Chen, Xinfu},
  title         = {Existence, uniqueness, and asymptotic stability of traveling waves in nonlocal evolution equations},
  journal       = {Adv. Differential Equations},
  volume        = {2},
  number        = {1},
  year          = {1997},
  pages         = {125--160},
  doi           = {10.57262/ade/1366809230},
}

@article{AK15,
  author        = {Achleitner, Franz and Kuehn, Christian},
  title         = {Traveling waves for a bistable equation with nonlocal diffusion},
  journal       = {Adv. Differential Equations},
  volume        = {20},
  number        = {9--10},
  year          = {2015},
  pages         = {887--936},
  doi           = {10.57262/ade/1435064517},
  url           = {http://projecteuclid.org/euclid.ade/1435064517},
}

@article{MR24,
  author        = {Mielke, Alexander and Reichelt, Sina},
  title         = {Traveling fronts in a reaction--diffusion equation with a memory term},
  journal       = {J. Dynam. Differential Equations},
  volume        = {36},
  number        = {S1},
  year          = {2024},
  pages         = {S487--S513},
  doi           = {10.1007/s10884-022-10133-6},
}

@article{BN91,
  author        = {Berestycki, Henri and Nirenberg, L.},
  title         = {On the method of moving planes and the sliding method},
  journal       = {Bol. Soc. Brasil. Mat. (N.S.)},
  volume        = {22},
  number        = {1},
  year          = {1991},
  pages         = {1--37},
  doi           = {10.1007/BF01244896},
}

\end{document}